\documentclass[11pt]{article}
\usepackage{amssymb}
\usepackage{amsmath}
\usepackage{latexsym}
\usepackage{amsthm}
\usepackage{mathptmx}

\makeatletter
\def\thmhead@plain#1#2#3{%
  \thmname{#1}\thmnumber{\@ifnotempty{#1}{ }#2}%
  \thmnote{ \the\thm@notefont(#3)}}
\let\thmhead\thmhead@plain
\def\swappedhead#1#2#3{%
  \thmnumber{#2}\thmname{\@ifnotempty{#2}{. }#1}%
  \thmnote{ \the\thm@notefont(#3)}}
\makeatother

\theoremstyle{definition} %%% for statements in roman typeface

 \newtheorem{definition}{Definition}[section]
 \newtheorem{remark}[definition]{Remark}

\theoremstyle{plain}      %%% for statements in italic typeface

 \newtheorem{proposition}[definition]{Proposition}
 \newtheorem{statement}[definition]{Statement}
 \newtheorem{theorem}[definition]{Theorem}
 
 \newtheorem{lemma}[definition]{Lemma}
 \newtheorem{defi}[definition]{Definition}
 
 \def \dem{\noindent{\sc Proof.~}}
\def \findem {\hfill{\hbox {\vrule\vbox{\hrule width 6pt\vskip 6pt\hrule}\vrule}}}

\def\OC{{\mathcal{O}}}
\def\QC{{\mathcal{Q}}}
\def\RC{{\mathcal{R}}}

\def\PC{{\mathcal{P}}}
\def\TC{{\mathcal{T}}}

\def\g{{\mathfrak{g}}}

\def\m{{\mathfrak{m}}}
\def\tf{{\mathfrak{t}}}
\def\Ff{{\mathfrak{F}}}

\def\fgam{{\mathfrak{S}}}

\def\NM{{\mathbb{N}}}
\def\QM{{\mathbb{Q}}}
\def\KM{{\mathbb{K}}}

\def \Bbo{\mathop{{\bigotimes}}\limits}
\def \ho{\widehat{\otimes}}
\def \bo{\bar{\otimes}}

\def \uk{\underline{k}}

\def \db {{\mathbf{\mathop{d}}}}
\def \T {\mathop{\hbox{\rm T}}\nolimits}

\def \gr {\mathop{\hbox{\rm gr}}\nolimits}

\def \Card {\mathop{\hbox{\rm card}}\nolimits}
\def \Ker {\mathop{\hbox{\rm Ker}}\nolimits}
\def \Hom {\mathop{\hbox{\rm Hom}}\nolimits}

\def \Pent {\mathop{\underline{\hbox{\rm Pent}}}\nolimits}
\def \Assocb {\mathop{\underline{\hbox{\rm Assoc}}}\nolimits}

\def \Aut {\mathop{\hbox{\rm Aut}}\nolimits}

\def \Alt {\mathop{\hbox{\rm Alt}}\nolimits}
\def \Ad {\mathop{\hbox{\rm Ad}}\nolimits}
\def \Im {\mathop{\hbox{\rm Im}}\nolimits}

\def \mod {\mathop{\hbox{\rm mod}}\nolimits}
\def \ad {\mathop{\hbox{\rm ad}}\nolimits}

\def \ptree {{\scriptstyle{\rm \hbox{-}tree}}}
\def \pa {{\scriptstyle{\rm \ a\ }}}
\def \pecard {{\scriptstyle{\rm card}}}

\def \dif {{\scriptstyle{\rm d}}}

\def \Lie {{\scriptstyle{\rm Lie}}}

\def \WX {{\scriptstyle{\rm WX}}}

\def \lef {{\scriptstyle{\rm left}}}

\def \adm {{\scriptstyle{\rm adm}}}
\def \Id {\mathop{\hbox{\rm id}}\nolimits}
\def \Oplus {\mathop{\oplus}\limits}

\def \CAP {\mathop{\cap}\limits}
\def \hopl {\mathop{\widehat{\oplus}}\limits}
\def \To {\mathop{\longrightarrow}\limits}
\def \Tto {\mathop{\longleftarrow}\limits}

\def \Lim {\mathop{\lim}\limits}

\def \Otimes {\mathop{\otimes}}

\def \T {\mathop{\hbox{\rm Tree}}\nolimits}

\def \dif{{\rm d}}

\begin{document}
\title{Poisson algebras associated to quasi-Hopf algebras}

\author{Benjamin Enriquez
and Gilles Halbout\\%}\address{
{\small{}}\cr
{\small{ \! Institut de Recherche Math\'ematique Avanc\'ee de Strasbourg}}\cr
{\small{UMR 7501 de l'Universit\'e Louis Pasteur et du CNRS}}\cr
{\small{7, rue R. Descartes F-67084 Strasbourg }}\cr
{\small{ \! e-mail:\,\texttt{enriquez@math.u-strasbg.fr}}}\cr
{\small{ \! e-mail:\,\texttt{halbout@math.u-strasbg.fr}}}\cr
}
\markboth{Benjamin Enriquez and Gilles Halbout}
{Poisson algebras associated to quasi-Hopf algebras}

\maketitle

\abstract{{\small We define admissible quasi-Hopf quantized universal 
enveloping (QHQUE) algebras by $\hbar$-adic valuation conditions. 
We show that any QHQUE algebra is twist-equivalent to an admissible one.
We prove a related statement: any associator is twist-equivalent 
to a Lie associator. We attach a quantized formal series algebra 
to each admissible QHQUE algebra and study the resulting Poisson algebras.
}}

\vskip20pt

\centerline {\bf \S\; 0 \ Introduction}

\vskip20pt

In \cite{WX}, Weinstein and Xu introduced a geometric counterpart of
quasitriangular quantum groups: they proved that if $(\g,r)$ is a
finite dimensional quasi-triangular Lie bialgebra, then the dual group
$G^*$ is equipped with a braiding $\RC_\WX$ with properties analogous
to those of quantum $R$-matrices (in particular, it is a set-theoretic
solution of the quantum Yang-Baxter Equation). An explicit relation to
the theory of quantum groups was later given in \cite{GH,EH,EGH}: to a
quasi-triangular QUE algebra $(U_\hbar(\g),m,R)$ quantizing $(\g,r)$,
one associates its quantized formal series algebra (QFSA)
$U_\hbar(\g)' \subset U_\hbar(\g)$; $U_\hbar(\g)'$ is a flat
deformation of the Hopf-co-Poisson algebra $\OC_{G^*}=(U(\g^*))^*$ of
formal functions of $G^*$. Then one proves that $\Ad(R)$ preserves
$U_\hbar(\g)'^{\bo 2}$, and $\Ad(R)|_{\hbar=0}$ coincides with the
automorphism $\RC_\WX$ of $\OC_{G^*}^{\bo 2}$; moreover, $\rho=\hbar
\log(R)|_{\hbar=0}$ is a function of $\OC_{G^*}^{\bo 2}$, independent
on a quantization of $\g^*$, which may be expressed universally in
terms of $r$, and $\RC_\WX$ coincides with the ``time one
automorphism'' of the Hamiltonian vector field generated by $\rho$.

\smallskip

In this paper, we study the analogous problem in the case of
quasi-quantum groups (quasi-Hopf QUE algebras).  The classical limit
of a QHQUE algebra is a Lie quasi-bialgebra (LQBA).  V. Drinfeld
proposed to attach Poisson-Lie ``quasi-groups'' to each LQBA
(\cite{Dr4}).  Axioms for Poisson-Lie quasi-groups are the quasi-Hopf
analogues of the Weinstein-Xu axioms.

A {\it Poisson-Lie quasi-group} is a Poisson manifold $X$, together
with a ``product'' Poisson map $X^2 \To^{m_X} X$, a unit for this
product $e\in X$, and Poisson automorphisms $\Phi_X \in \Aut(X^3)$,
$\Phi_X^{12,3,4}$, $\Phi_X^{1,23,4}$ and $\Phi_X^{1,2,34}\in
\Aut(X^4)$, such that
\begin{align*}
m_X  \circ (\Id \times  m_X)&=m_X \circ (m_X \times \Id) \circ \Phi_X, \\
(m_X \times \Id\times \Id) \circ \Phi_X^{12,3,4}&=\Phi_X \circ (m_X \times \Id
\times \Id),\\
(\Id \times m_X \times \Id) \circ \Phi_X^{1,23,4}&=\Phi_X \circ
(\Id \times m_X \times \Id),~\hbox{etc.}\\
\hbox{and }~ \Phi_X^{1,2,34} \circ \Phi_X^{12,3,4}&=
(\Id \times \Phi_X) \circ \Phi_X^{1,23,4} \circ (\Phi_X \times \Id).
\end{align*}
A {\it twistor} for the quasi-group $(X,m_X,\Phi_X)$ is a collection of
Poisson 
automorphisms $F_X \in \Aut(X^2)$, $F_X^{12,3}$, $F_X^{1,23} \in \Aut(X^3)$,
$F_X^{(12)3,4}$, $F_X^{1(23),4}$, $F_X^{12,34}$, $F_X^{1(23),4}$,
$F_X^{1,(23)4}\in \Aut(X^4)$ such that 
\begin{align*}
(m_X \times \Id) \circ F_X^{12,3}&=F_X \circ (m_X \times \Id),\\
\big( (m_X \circ (\Id \times m_X)) \times \Id \big) \circ F_X^{1(23),4}&=
F_X \circ \big( (m_X \circ (\Id \times m_X))\times \Id \big),\\
F_X^{(12)3,4}&=(\Phi_X \times \Id) \circ F_X^{1(23),4} \circ
(\Phi_X \times \Id)^{-1},\hbox{ etc.}
\end{align*}
A twistor replaces the quasi-group $(X,m_X,\Phi_X)$ by
$ (X,m_X',\Phi_X')$ with  $m_X'=m_X \circ F_X$ and 
$\Phi_X'=(F_X^{1,23})^{-1}\circ (F_X \times \Id)^{-1} \circ \Phi_X \circ
F_X^{1,23} \circ (\Id \times F_X).$

(Other axioms for Poisson-Lie quasi-groups were proposed
in a differential-geometric language in \cite{Ban,KS}.)

\smallskip

We do not know a ``geometric'' construction of a twist-equivalence
class of $(X,m_X,\Phi_X)$ associated to each Lie quasi-bialgebra, in
the spirit of \cite{WX}.  Instead we generalize the ``construction of
a QFS algebra and passage to Poisson geometry'' part of the above
discussion, and we derive from there a construction of triples
$(X,m_X,\Phi_X)$, in the case of Lie quasi-bialgebras with vanishing 
cobracket. 

Let us describe the generalization of the ``construction of a QFS
algebra'' part (precise statements are in Section 1).  We introduce
the notion of an {\it admissible} quasi-Hopf QUE algebra, and we
associate a QFSA to such a QHQUE algebra. Each QHQUE algebra can be
made admissible after a suitable twist.

We generalize the ``passage to Poisson geometry'' part as follows.
The reduction modulo $\hbar$ of the obtained QFS algebra is a
quintuple $(A,m,P,\Delta,\widetilde{\varphi})$ satisfying certain
axioms; in particular $\exp(V_{\widetilde{\varphi}})$ is an
automorphism of $A^{\ho 3}$, and $(A,m,\exp(V_{\widetilde{\varphi}}))$
satisfies the axioms dual to those of $(X,m_X,\Phi_X)$.

When the Lie quasi-bialgebra arises from a metrized Lie algebra,
admissible QHQUE algebras quantizing it are given by Lie associators,
and we obtain a quasi-group $(X,m_X,\Phi_X)$ using our construction.
We also prove that its twist-equivalence class does not depend on the
choice of an associator.

Finally, we prove a related result: any associator is twist-equivalent
to a unique Lie associator.

\vskip20pt

\centerline {\bf \S\; 1 \ Outline of results}
\stepcounter{section}\label{section:1}

\vskip20pt

Let $\KM$ be a field of characteristic $0$.  Let $(U,m)$ be a
topologically free $\KM[[\hbar]]$-algebra equipped with algebra
morphisms
$$\Delta~:~U\to U\ho U,\hbox{ and }\epsilon~:~U \to \KM[[\hbar]]$$
$$\hbox{with }(\epsilon \otimes \Id)\circ \Delta=(\Id \otimes
\epsilon)\circ \Delta =\Id$$
such that the reduction of
$(U,m,\Delta) \hbox{ modulo } \hbar$ is a universal enveloping algebra.
Set
$$U'=\{x \in U|~\hbox{for any tree }P,~ \delta^{(P)}(x) \in
\hbar^{|P|}U^{\otimes |P|}\}$$
(see the definitions of a tree, $\delta^{(P)}$, and $|P|$ in Section 2).
We prove:
\begin{theorem}
\label{theo:1}
$U'$ is a topologically free $\KM [[\hbar]]$-algebra.
It is equipped with a complete decreasing algebra filtration
$$(U')^{(n)}=\{x \in U|~\hbox{for any tree }P,~
\delta^{(P)}(x) \in \hbar^n U^{\otimes |P|}\}.$$
$U'$ is stable under the multiplication $m$ and the map $\Delta$~: $U\to U^{\ho 2}$ induces
a continuous algebra morphism
$$\Delta_{U'}~:~U' \to {U'}^{\bo 2}=
\lim_{\Tto_n}\left({U'}^{\ho 2}/\sum_{p,q|p+q=n}{U'}^{(p)}\otimes
{U'}^{(q)}\right).$$
Set $\OC:=U'/\hbar U'$. Then $\OC$ is a complete commutative local ring and the
reduction modulo $\hbar$ of $\Delta_{U'}$ is a continuous ring morphism
$$\Delta_\OC~:~\OC \to \OC^{\bo 2}=
\lim\limits_{\Tto_n}\left(\OC^{\otimes 2}/\sum\limits_{p,q|p+q=n} \OC^{(p)}\otimes
\OC^{(q)}\right),$$ 
where $\OC^{(p)}={U'}^{(p)}/(\hbar U \cap {U'}^{(p)})$.
\end{theorem}
\begin{theorem}
\label{theo:2}
Let $(U,m,\Delta,\Phi)$ be a quasi-Hopf QUE algebra.
Assume that 
\begin{equation}
\label{admi}
 \hbar \log(\Phi) \in (U')^{\bo 3} .
\end{equation}
Then there is a noncanonical isomorphism
of filtered algebras
$U'/\hbar U' \to \widehat{S}^\cdot(\g),$
where $\widehat{S}^\cdot(\g)$ is the formal
series completion of the symmetric algebra $S^\cdot(\g)$.
\end{theorem}
\noindent
When $(U,m,\Delta,\Phi)$ satisfies the hypothesis (\ref{admi}), we say that it is
{\it admissible}. In that case, we say that $U'$ is the quantized formal series algebra
(QFSA) corresponding to $(U,m,\Delta,\Phi)$.
\noindent
Let us recall the notion of a {\it twist} of a quasi-Hopf QUE
algebra 
$(U,m,\Delta,\Phi)$. This is an element $F\in \big( U^{\ho 2}\big)^\times$,
such that $(\epsilon \otimes \Id)(F)=(\Id \otimes \epsilon)(F)=1$.
It transforms
$(U,m,\Delta,\Phi)$ into the quasi-Hopf algebra $(U,m,{}^F\Delta,{}^F\Phi)$,
where
$${}^F\Delta=\Ad(F) \circ \Delta,\hbox{ and }{}^F\Phi=
(1 \otimes F)(\Id \otimes \Delta)(F)\Phi(\Delta \otimes \Id)(F)^{-1}(F\otimes 1)^{-1}.$$
\begin{theorem}
$\hbox{~}$
\label{theo:3}
\begin{description}
\item[1)]
Let $(U,m,\Delta,\Phi)$ be an admissible quasi-Hopf QUE algebra. Let us
say that a twist $F$ of $U$ is admissible if  
$\hbar \log(F) \in {U'}^{\bo 2}\!$.
Then the twisted quasi-Hopf algebra
$(U,m,{}^F\Delta,{}^F\Phi)$ is also admissible, and its QFSA coincides with
$U'$.
\item[2)] Let $(U,m,\Delta,\Phi)$ be an arbitrary quasi-Hopf QUE algebra. There
exists a twist $F_0$ of $U$ such that the twisted quasi-Hopf algebra 
$(U,m,{}^{F_0}\Delta,{}^{F_0}\Phi)$ is admissible.
\end{description}
\end{theorem}
\noindent Theorem \ref{theo:3} can be interpreted as follows.
Let $(U,m)$ be a formal deformation of a universal enveloping algebra.
The set of twists of $U$ is a subgroup 
$\TC$ of $(U^{\ho 2})^\times$.
Denote by $\QC$ the set of all quasi-Hopf structures on
$(U,m)$, and by
$\QC_\adm$ the subset of admissible structures.
If $\QC$ is nonempty, then $\QC_\adm$ is also nonempty, and all its
elements
give rise to the same subalgebra
$U' \subset U$ (Theorem \ref{theo:3}, 1)).
Using $U'$, we then define the subgroup $\TC_\adm \subset \TC$ of
admissible
twists. We have a natural action of $\TC$ on $\QC$, which restricts to
an action
of $\TC_\adm$ on $\QC_\adm$. Theorem \ref{theo:3} 2) says that the
natural map
$$\QC_\adm/\TC_\adm \to \QC/\TC$$
is surjective. Let us explain why it is not injective in general.  Any
QUE Hopf algebra $(U,m,\Delta)$ is admissible as a quasi-Hopf algebra.
If $u \in U^\times$ and $F=(u \otimes u)\Delta(u)^{-1}$, then
$(U,m,{}^F\Delta)$ is a Hopf algebra.  So $(U,m,\Delta)$ and
$(U,m,{}^F\Delta)$ are in the same class of $\QC/\TC$.  These are also
two elements of $\QC_\adm$; the corresponding QFS algebras are $U'$
and $\Ad(u)(U')$.  In general, these algebras do not coincide, so
$(U,m,\Delta)$ and $(U,m,{}^F\Delta)$ are not in the same class of
$\QC_\adm/\TC_\adm$.

\medskip
The following result is a refinement of Proposition 3.10 of 
\cite{Dr2}. Let $(\g,\mu,\varphi)$ be a pair of a Lie algebra 
$(\g,\mu)$ and 
$\varphi\in\wedge^3(\g)^\g$. Then $(\g,\delta = 0,\varphi)$ is a 
Lie bialgebra. 
\begin{proposition} \label{refinement}
There exists a series ${\cal E}(\varphi) \in U(\g)^{\otimes 3}[[\hbar]]$, 
expressed in terms of $(\mu,\varphi)$ by universal acyclic expressions, 
such that $(U(\g)[[\hbar]],m_0,\Delta_0,{\cal E}(\varphi))$ is an 
admissible quantization of $(\g,\mu,\varphi)$. 
\end{proposition}
This proposition is proved in Section 6. 
\medskip 

Let us define a Drinfeld algebra as follows:
\begin{defi}
A Drinfeld algebra is a quintuple
$(A,m_0,P,\Delta,\widetilde{\varphi})$,
where 
\begin{itemize}
\item $(A,m_0)$ is a formal series algebra, 
\item $P$ is a Poisson strcture
on $A$ ``vanishing at the origin''
(i.e., such that $\Im(P) \subset \m_A$, where
$\m_A$ is the maximal ideal of $A$),
\item
$\Delta$~:~$A \to A \ho A$ is a continuous
Poisson algebra morphism, such that
$(\epsilon \otimes \Id)\circ \Delta = (\Id \otimes \epsilon) \circ
\Delta = \Id$, where $\epsilon$~: $A \to A/\m_A=\KM$ is the natural
projection, 
\item $\widetilde{\varphi}\in (\m_A)^{\ho 3}$ satisfies
\begin{align*}
(\Id \otimes \Delta)(\Delta(a))&=\widetilde{\varphi}
\star (\Delta \otimes \Id)(\Delta(a))\star (-\widetilde{\varphi}),~a
\in A,\\
\widetilde{\varphi}^{1,2,34} \star\
\widetilde{\varphi}^{12,3,4}&=
\widetilde{\varphi}^{2,3,4}
\star \ \widetilde{\varphi}^{1,23,4}
\star\ \widetilde{\varphi}^{1,2,3},
\end{align*}
where we set $f \star g = f + g + {{1}\over{2}}P(f,g) + \cdots$, the
Cambell-Baker-Hausdorff (CBH) series of the Lie algebra $(A,P)$.
\end{itemize}
\end{defi}
\noindent
If $\widetilde{f} \in \m_A^{\ho 2}$, we define the twist of 
the Drinfeld algebra $(A,m_0,P,\Delta,\widetilde{\varphi})$ by $\widetilde{f}$
as the algebra
$(A,m_0,P,{}^{\widetilde{f}}\Delta,{}^{\widetilde{f}}\widetilde{\varphi})$,
where
\begin{align*}
{}^{\widetilde{f}}\Delta(a)&=\widetilde{f}\star \Delta(a)
\star(-\widetilde{f}),\hbox{ and}\\
{}^{\widetilde{f}}\widetilde{\varphi}&=\widetilde{f}^{2,3}\star
\widetilde{f}^{1,23}\star \widetilde{\varphi}\star
(-\widetilde{f}^{12,3})\star
(-\widetilde{f}^{1,2});
\end{align*}
then
$(A,m_0,P,{}^{\widetilde{f}}\Delta,{}^{\widetilde{f}}\widetilde{\varphi})$
is again a Drinfeld algebra.

\smallskip

\begin{remark}
If $\Lambda$ is any Artinian local $\KM$-ring with residue field
$\KM$, set $X=\Hom_\KM(A,\Lambda)$.  Then $X$ is the ``Poisson-Lie
quasi-group'', in the sense of the Introduction.  Namely, $\Delta_0$
induces a product $m_X$~: $X\times X \to X$, and
$\exp(V_{\widetilde{\varphi}})$,
$\exp(V_{\widetilde{\varphi}^{12,3,4}})$, etc., induce automorphisms
$\Phi_X$, $\Phi_X^{12,3,4}$, etc., of $X$, that satisfy the
quasi-group axioms (we denote by $V_f$ the Hamiltonian derivation of
$A^{\ho k}$ induced by $f \in A^{\ho k}$). Moreover, if
$\widetilde{f}$ is a twist of $A$, then $\exp(V_{\widetilde{f}})$,
$\exp(V_{\widetilde{f}^{12,3}})$, $\exp(V_{\widetilde{f}^{(12)3,4}})$,
etc., define a twistor $(F_X,F_X^{12,3},F_X^{(12)3,4},\dots)$ of
$(X,m_X,\Phi_X)$.  Twisting $A$ by $\widetilde{f}$ corresponds to
twisting $(X,m_X,\Phi_X)$ by $(F_X,F_X^{12,3},\dots)$.
\end{remark}
\begin{lemma}
\label{reduction:lemma}
If $(A,m_0,P,\Delta,\widetilde{\varphi})$ is a Drinfeld algebra, set
$\g=\m_A/(\m_A)^2$; then $P$ induces a Lie bracket $\mu$ on $\g$,
$\Delta- \Delta^{1,2}$ induces a linear
map $\delta$~: $\g \to \Lambda^2(\g)$, and
the reduction of $\Alt(\widetilde{\varphi})$ is an element $\varphi $
of $\Lambda^3(\g)$.
Then $(\g,\mu,\delta,\varphi)$ is a Lie quasi-bialgebra.\ Moreover, twisting
$(A,m_0,P,\Delta,\widetilde{\varphi})$ by $\widetilde{f}$ corresponds
to twisting $(\g,\mu,\delta,\varphi)$ by
$$f := \big(\Alt(\widetilde{f}) \mod \ (\m_A)^2 \otimes \m_A +
\m_A \otimes (\m_A)^2\big)\in \Lambda^2(\g).$$
\end{lemma}
\noindent
Taking the reduction modulo $\hbar$ induces a natural map
$$\QC_\adm/\TC_\adm \to \{\hbox{Drinfeld algebra structures on }
\widehat{S}^\cdot(\g)\}/\hbox{twists}.$$
To summarize, we have a diagram
$$\begin{matrix}
\QC/\TC &\leftarrow& \QC_\adm/\TC_\adm &\to& \left\{\begin{matrix}
%\scriptstyle{\rm{Drinfeld~ algebra ~ structures}}\\
%\scriptstyle{\rm{on~ }}\widehat{S}^\cdot(\g)
\hbox{{Drinfeld~algebra~structures}}\\
\hbox{{on~ }}\widehat{S}^\cdot(\g)
\end{matrix}\right\}/\hbox{twists}\cr
{\scriptstyle{\it{class~}}}\downarrow&&&&\downarrow{\scriptstyle{\it{~red}}}
\end{matrix}$$
$$\{\hbox{Lie quasi-bialgebra structures on }
(\g,\mu)\}/\hbox{twists},\hskip2cm$$
where {\it class} is the classical limit map described in \cite{Dr2}, and
{\it red} is the map described in Lemma \ref{reduction:lemma}. It is easy to
see
that this diagram commutes.

\medskip

When $U$ is a Hopf QUE algebra, the corresponding Drinfeld
algebra is the Hopf-Poisson structure on
$\OC_{G^*}=(U(\g^*))^*$, and
$\widetilde{\varphi}=0$.

\medskip

Let $(\g,\mu,\delta,\varphi)$ be a Lie quasi-bialgebra. A 
{\it lift of }$(\g,\mu,\delta,\varphi)$ is a Drinfeld algebra, whose
reduction is $(\g,\mu,\delta,\varphi)$.
A general problem is to  construct a lift for any Lie
quasi-bialgebra.
We will not solve this problem, but we will give partial existence and
unicity results. 

Assume that $\delta = 0$. A quasi-Lie bialgebra is then the same as 
a triple $(\g,\mu,\varphi)$ of a Lie algebra $(\g,\mu)$
and $\varphi\in \wedge^3(\g)^\g$. 

\begin{theorem}\label{theo:4}
$\hbox{~}$
\begin{description}
\item[1)] There exists a lift 
\begin{equation} \label{num}
(\widehat{S}^\cdot(\g),m_0,P_{\g^*},\Delta_0,\widetilde\varphi) 
\end{equation}
of $(\g,\mu,\delta = 0,\varphi)$. Here $P_{\g^*}$ is
the Kostant-Kirillov Poisson structure on $\g^*$ and $\Delta_0$ is the
coproduct for which the elements of $\g$ are primitive.

\item[2)] Any two lifts of $(\g,\mu,\delta=0,\varphi)$ of the form 
(\ref{num}) are related by a $\g$-invariant twist. 
\end{description}
\end{theorem}

Examples of quasi-Lie bialgebras with $\delta = 0$ arise from 
{\it metrized Lie algebras}, i.e., pairs $(\g,t_\g)$ of
a Lie algebra $\g$ and $t_\g \in S^2(\g)^\g$.
Then $\varphi=[t_\g^{1,2},t_\g^{2,3}]$.
Recall that a {\it Lie associator} is a noncommutative formal series
$\Phi(A,B)$, such that $\log\Phi(A,B)$ is a Lie series $[A,B]+$higher
degrees terms, satisfying the pentagon and hexagon identities 
(see \cite{Dr3}).

\begin{proposition} \label{metrized}
If $\Phi$ is a Lie associator, we may set $\varphi =
\log(\Phi)(\bar{t}_\g^{1,2},\bar{t}_\g^{2,3})$, where 
$\bar{t}_\g^{i,j}$ is the image of $t_\g^{i,j}$ in
$\widehat{S}^\cdot(\g)^{\ho 3}$, and we use the Poisson bracket of
$\widehat{S}^\cdot(\g)^{\ho 3}$ in the expression
of $\log (\Phi)(\bar{t}_\g^{1,2},\bar{t}_\g^{2,3})$.
\end{proposition}
\noindent
We prove these results in Section 6. If now $\Phi$ is a general
(non-Lie) associator, $(U(\g)[[\hbar]],$
$m_0,\Delta_0,\Phi(\hbar t_\g^{1,2},\hbar t_\g^{2,3}))$ is a quasi-Hopf QUE algebra, but it is
admissible only when
$\Phi$ is Lie (for general $\g$). According to Theorem \ref{theo:3} 2), it is
twist-equivalent to an admissible quasi-Hopf QUE algebra.
We prove
\begin{theorem}
\label{theo:7.1}
Any (non-Lie) associator is twist-equivalent
to a unique Lie associator.
\end{theorem}
\noindent So the ``concrete'' version of the twist of Theorem \ref{theo:7.1}
is an example of the twist $F$ of Theorem \ref{theo:3} 2).

\vskip20pt

\centerline {\bf \S\; 2 \ Definition and properties of $U'$}

\stepcounter{section}\label{section:2}

\vskip20pt

In this section, we prove Theorem \ref{theo:1}. We first introduce the
material for the
definition of $U'$: trees (a); the map $\delta^{(P)}$ (b); then
we prove Theorem \ref{theo:1} in (c) and (d).

\bigskip

\noindent {\bf - a - Binary complete planar rooted trees}

\bigskip

\begin{defi}
A $n$-binary complete planar rooted tree ($n$-tree for short) is a set
of vertices and oriented edges satisfying the following conditions:
\begin{itemize}
\item each edge carries one of the labels $\{l,r\}$.
\item if we set: 
$$\hbox{valency of a vertex}=(\Card(\hbox{incoming edges}),\Card(\hbox{outgoing
edges})),$$ 
we have
\begin{itemize}
\item there exists exactly one vertex with valency $(0,2)$ (the root)
\item there exists exactly $n$ vertices with valency $(1,0)$ (the leaves)
\item all other vertices have valency $(1,2)$
\item if a vertex has valency $(x,2)$, then one of its outgoing edges has label
$l$ and the other has label $r$.
\end{itemize}
\item the set of leaves has cardinal $n$. 
\end{itemize}
\end{defi}
\noindent
Let us denote, for $n \geq 2$, 
$$\T_n=\{n\hbox{-binary complete planar rooted
trees}\}.$$
By definition, $\T_1$ consists of one element (the tree with a root
and one nonmarked edge) and $\T_0$ consists of one element
(the tree with a root and no edge). We will write $|P|=n$ if $P$ is a tree in $\T_n$.
\begin{defi}
(Extracted trees) Let $P$ be a binary complete planar rooted tree.
Let $L$ be the set of its leaves and let $L'$ be a subset of $L$. 
We define the extracted subtree $P_{L'}$ as follows:
\begin{description} 
\item[(1)] $\widetilde{P}_{L'}$ is the set of all edges connecting the root with an
element of $L'$,
\item[(2)] the vertices of $\widetilde{P}_{L'}$ all have valency $(0,2)$, $(1,0)$,
$(1,2)$ or $(1,1)$;
\item[(3)] $P_{L'}$ is obtained from $\widetilde{P}_{L'}$ by replacing each maximal sequence
of edges related by a $(1,1)$ vertex, by a single edge whose label is the label of the
first edge of the sequence.
\end{description}
Then $P_{L'}$ is a $|L'|$-binary complete planar rooted tree.
\end{defi}
\begin{defi}
(Descendants of a tree)
If we cut the tree $P$ by removing its root and the related vertices, we get two
trees $P'$ and $P''$, its {\em left and right descendants}.
\end{defi}

In the same way, we define the left and right descendants of a vertex of $P$. 

If $P$ is a $n$-tree, there exists a unique bijection of the set of leaves
with $\{1,\dots,n\}$, such that for each vertex, the number attached to any leaf of
its left descendant is smaller than the number attached to any leaf of its right descendant.

\bigskip

\noindent {\bf - b - Definition of $\delta^{(P)}~: U \to U^{\ho n}$}

\bigskip

Let us place ourselves in the hypothesis of Theorem \ref{theo:1}. 
Let us define $\delta^{(2)}$~: $U \to U^{\ho 2}$,
$$\delta^{(2)}(x)=\Delta(x)-x\otimes 1 -1 \otimes x + \epsilon(x) 1 \otimes 1.$$
For $P_2$ the only tree of $\T_2$, we set 
$$\delta^{(P_2)}=\delta^{(2)}=\delta.$$
For $P_1$, the only tree of $\T_1$, we set 
$$\Delta^{(P_1)}(x)=\delta^{(1)}(x)=x
-\epsilon(x)1.$$
For $P_0$ the only tree of $\T_0$,
we set 
$$\delta^{(P_0)}(x)=\delta^{(0)}(x)=\epsilon(x).$$
When $P$ is a $n$-tree with descendants $P'$ and $P''$, we set
$$
\delta^{(P)}=(\delta^{(P')} \otimes \delta^{(P'')})\circ \delta,
$$
so $\delta^{(P)}\hbox{ is a linear map }U \to U^{\ho n}$.

\bigskip

\noindent {\bf - c - Behavior of $\delta^{(P)}$ with respect to multiplication}

\bigskip

If $\Sigma=\{i_1,\dots,i_k\}$ is a subset of $\{1,\dots,n\}$, where $i_1 <i_2 <
\cdots <i_k$, the map $x \mapsto x^\Sigma$ is the linear map
$U^{\ho k} \to U^{\ho n}$, defined by
$$x_1\otimes \cdots \otimes x_k \mapsto 1^{\otimes i_1-1} \otimes x_1 \otimes 
1^{\otimes i_2-i_1-1} \otimes x_2 \otimes \cdots \otimes 1^{\otimes i_k -
i_{k-1}-1} \otimes x_k \otimes 1^{\otimes n-i_k-1}.$$
If $\Sigma=\emptyset$, $x \mapsto x^\Sigma$ is the map
$\KM \to U^{\ho n}$, $1 \mapsto 1^{\otimes n}$.

\begin{proposition}
\label{prop:1.4}
\label{multiplication}
For $P \in \T_n$, we have the identity
$$\delta^{(P)}(xy)=\sum_{
{\begin{array}{l}
{\scriptstyle{\Sigma',\Sigma''\subset \{1,\dots,n\}|}}\cr
{\scriptstyle{\Sigma'\cup\Sigma''=\{1,\dots, n\}}}
\end{array}}
}(\delta^{({\Sigma'})}(x))^{\Sigma'}(\delta^{({\Sigma''})}(y))^{\Sigma''},$$
for any $x,y \in U$.
\end{proposition}
\noindent This proposition is proved in Section 5.

\bigskip

\noindent {\bf - d - Construction of $U'$}

\bigskip

Let us set
$$U'=\{x \in U|~\hbox{for any tree }P,~ \delta^{(P)}(x) \in
\hbar^{|P|}U^{\ho |P|}\}.$$
Then $U'$ is a topologically free $\KM[[\hbar]]$-submodule of $U$.
Moreover, if $x,y \in U'$, and
$P$ is a tree, then
$$\delta^{(P)}([x,y])=
\sum_{\begin{array}{l}
{\scriptstyle{\Sigma,\Sigma' \subset \{1,\dots,|P|\}}}\cr
{\scriptstyle{\Sigma \cup\Sigma' = \{1,\dots,|P|\}}}
\end{array}}
\left[\delta^{(P_{\Sigma})}(x)^{\Sigma},
\delta^{(P_{\Sigma'})}(y)^{\Sigma'}\right];$$
the summand corresponding to a pair $(\Sigma,\Sigma')$ with $\Sigma \cap
\Sigma'=\emptyset$ is zero, and the $\hbar$-adic valuation of
the other summands is $\leq |\Sigma| + |\Sigma'| \leq |P| +1$; so
$\delta^{(P)}([x,y]) \in \hbar^{|P|+1}U^{\ho |P|}$.
On the other hand, there exists $z \in U$ such
that 
$[x,y]=\hbar z$, so 
$\delta^{(P)} (z) \in \hbar^{|P|}U^{\ho |P|}$; so
$z \in U'$ and
we get $[x,y] \in \hbar U'$.
It follows that $U'/\hbar U'$ is commutative. Let us set
\begin{equation}
\label{filtr:1}
{U'}^{(n)}=U'\cap \hbar^nU.
\end{equation}
We have a decreasing filtration
$$U' = {U'}^{(0)} \supset 
{U'}^{(1)} \supset 
{U'}^{(2)} \supset \cdots ;$$
we have ${U'}^{(n)} \subset \hbar^n U$, so $U'$ is complete for the 
topology induced by this filtration. This
is an algebra filtration, 
i.e., ${U'}^{(i)}{U'}^{(j)}\subset
{U'}^{(i+j)}$. It induces an algebra filtration on
$U'/\hbar U'$,
$$U'/\hbar U' \supset \cdots \supset {U'}^{(i)}/\big(
{U'}^{(i)}\cap \hbar U'\big) \supset \cdots ,$$
for which $U'/\hbar U'$ is complete. 
Moreover, the completed tensor product
$$U' \bo U'=\lim_{\Tto_n} \big( U' \ho U'/\sum_{p,q|p+q=n}{U'}^{(p)}\ho
{U'}^{(q)}\big)$$
identifies with
\begin{multline*}
\lim_{\Tto_n}\big(
\{x \in U \ho U| \forall P, Q,(\delta^{(P)} \otimes \delta^{(Q)})(x)
\in \hbar^{|P|+|Q|}U^{\ho 2}\}/ \\
\{x \in U\ho U| \forall P,Q, ~(\delta^{(P)} \otimes \delta^{(Q)})(x)
\in \hbar^{\max(n,|P|+|Q|)}U^{\ho 2}\}\big).
\end{multline*}
If $x \in U'$, and $P,Q$ are trees, with $|P|,|Q|\not= 0$,
then since $\delta^{(P)}(1)=\delta^{(Q)}(1)=0$, we have
\begin{multline*}
(\delta^{(P)}\otimes \delta^{(Q)})(\Delta(x))=
(\delta^{(P)}\otimes \delta^{(Q)})(\delta(x))=
\delta^{(R)}(x) \in \hbar^{|R|}U^{\ho |R|}\\
=\hbar^{|P|+|Q|}U^{\ho |P|+|Q|},
\end{multline*}
where $R$ is the tree whose left and right descendants are $P$ and $Q$;
so $|R|=|P|+|Q|$.
On the other hand,
\begin{align*}
(\delta^{(P)}\otimes \epsilon)(\Delta(x))&=\delta^{(P)}(x)\otimes 1
\in \hbar^{|P|}U^{\otimes |P|}\\
(\epsilon \otimes \delta^{(P)})(\Delta(x))&=1 \otimes \delta^{(P)}(x)
\in \hbar^{|P|}U^{\otimes |P|},
\end{align*}
so $\Delta(x)$ satisfies 
$(\delta^{(P)}\otimes \delta^{(Q)})(\Delta(x))
\in \hbar^{|P|+|Q|}U^{\bo |P|+|Q|}$ for any pair of trees $(P,Q)$.
$\Delta$~: $U \to U \ho U$ therefore induces an algebra morphism
$\Delta_{U'}$~: $U' \to {U'}^{\bo 2}$, whose reduction modulo $\hbar$
is a morphism of complete local rings
$$\OC \to \OC^{\bo 2}=\lim_{\Tto_n}\left(\OC^{\otimes 2}/\sum_{p,q|p+q=n}\OC_p
\otimes
\OC_q\right),$$
where $\OC=U'/\hbar U'$ and $\OC_p={U'}^{(p)}/({U'}^{(p)}\cap \hbar U').$
\vskip20pt

\centerline {\bf \S\; 3 \ Classical limit of $U'$}

\stepcounter{section}\label{section:3}

\vskip20pt

We will prove Theorem \ref{theo:2}
as follows. We first compare the various $\delta^{(P)}$,
where $P$ is a $n$-tree (Proposition \ref{rels:delta}).
Relations found between the $\delta^{(P)}$ imply that they have $\hbar$-adic valuation properties
close to those of the Hopf case (Proposition \ref{prepa:flatness}).
We then prove Theorem \ref{theo:2}.
\bigskip

\noindent {\bf - a - Comparison of the various $\delta^{(P)}$}

\bigskip

Let $P$ and $P_0$ be $n$-trees. There exists an element
$\Phi^{P,P_0} \in U^{\ho n}$, such that
$\Delta^{(P)}=\Ad(\Phi^{P,P_0})\circ \Delta^{(P_0)}$.
The element $\Phi^{P,P_0}$ is a product of images of $\Phi$ and $\Phi^{-1}$
by the various maps
$U^{\ho 3} \to U^{\ho n}$ obtained by iteration of $\Delta$.
We have 
\begin{equation}
\label{transitivity}
\Phi^{P',P_0}=\Phi^{P',P}\Phi^{P,P_0}
\end{equation}
for any $n$-trees $P_0,P,P'$. For example,
\begin{align*}
(\Id \otimes \Delta)\circ \Delta &= \Ad(\Phi) \circ ((\Delta \otimes \Id)\circ
\Delta),\\
(\Delta \otimes \Delta)\circ \Delta &= \Ad(\Phi^{12,3,4}) \circ ((\Delta \otimes 
\Id^{\otimes 2})\circ (\Delta \otimes \Id)\circ
\Delta), \hbox{ etc.}
\end{align*}
\begin{proposition}
\label{rels:delta}
Assume that $\hbar \log(\Phi)\in (U')^{\bo 3}$.
Then there exists a sequence of elements
$$F^{PP_0R\Sigma \nu}=
\sum_\alpha F_{1,\alpha}^{PP_0R\Sigma\nu}\otimes \cdots
\otimes
F_{\nu,\alpha}^{PP_0R\Sigma\nu}\in ({U'}^{\bo n})^{\bo \nu},$$
indexed by the triples $(R,\Sigma,\nu)$, where $R$ is a tree such that $|R|<n$,
$\Sigma$ is a subset of $\{1,\dots,n\}$ with $\Card(\Sigma)=|R|$,
and $\nu$ is an integer 
$\geq 1$, such that the equality
\begin{multline}
\label{delta:P:P0}
\delta^{(P)}=\Ad(\Phi^{P,P_0})\circ \delta^{(P_0)}
+\sum_{k|k<n}\sum_{R \pa k\ptree}\sum_{\begin{array}{l}
{\scriptstyle{\Sigma \subset\{1,\dots,n\},}}\cr
{\scriptstyle{\pecard(\Sigma)=k}}
\end{array}}\\
\sum_{\nu \geq 1} \sum_{\alpha}\ad_\hbar(F_{1,\alpha}^{PP_0R\Sigma\nu})\circ \cdots \circ \ad_\hbar
(F_{\nu,\alpha}^{PP_0R\Sigma\nu})\circ (\delta^{(R)})^\Sigma
\end{multline}
holds. Here $\ad_\hbar(x)(y)={{1}\over{\hbar}}[x,y]$.
\end{proposition}
\dem
Let us prove this statement by induction on $n$. When $n=3$, we find
$$\delta^{(1(23))}=\Ad(\Phi)
\delta^{((12)3)}+(\Ad(\Phi)-1)
(\delta^{1,2}+\delta^{1,3}+\delta^{2,3}+\delta^{(1)1}+\delta^{(1)2}+\delta^{(1)3}),$$
so the identity holds with $F^{PP_0R\Sigma\nu}={{1}\over{\nu !}}
(\hbar\log \Phi)^{\bo \nu}$ for all choices
of $(R,\Sigma,\nu)$, except when $|R|=0$, in which case $F^{PP_0R\Sigma\nu}=0$.
Assume that the statement holds for any pair of $k$-trees, $k \leq n$, and let
us prove it for a pair $(P,P_0)$ of $(n+1)$-trees.
For $k$ any integer, let $P_\lef(k)$ be the $k$-tree corresponding to
$$\delta^{(P_\lef(k))}=(\delta \otimes \Id^{\otimes k-2})\circ
\cdots \circ  \delta.$$
Thanks to (\ref{transitivity}), we may assume that $P_0=P_\lef(n+1)$ and $P$ is
arbitrary. Let $P'$ and $P''$ be the subtrees of $P$, such that
$|P'|+|P''|=n+1$,
and $\delta^{(P)}=(\delta^{(P')}\otimes \delta^{(P'')})\circ \delta$. Let $P_1$
and $P_2$ the $n$-trees such that
$$\delta^{(P_1)}=(\delta^{(P_\lef(k'))}\otimes \delta^{(P'')})\circ \delta\hbox{
and }
\delta^{(P_2)}=(\delta^{(P_\lef(k'))}
\otimes \delta^{(P_\lef(k''))})\circ \delta$$
Assume that $|P_1|\not= 1$. Using (\ref{transitivity}), we reduce the proof of
(\ref{delta:P:P0}) to the case of the pairs $(P,P_1),(P_1,P_2)$ and
$(P_2,P_0)$.
Then the induction hypothesis applied to the pair $(P'\!,P_\lef(k'))$,
together with $\Phi^{P,P_1}=\Phi^{P',P_\lef(k')}\otimes 1^{\otimes k''}$,
implies
\begin{multline*}
\delta^{(P)}=\Ad(\Phi^{P,P_1})\circ \delta^{(P_1)}+ \sum_{k|k<k'}
\sum_{R\pa k\ptree}\sum_{\begin{array}{l}
{\scriptstyle{\Sigma \subset\{1,\dots,k'\},}}\cr
{\scriptstyle{\pecard(\Sigma)=k}}
\end{array}}\\
\sum_{\nu \geq 1} \sum_\alpha
\Ad(\Phi^{P,P_1}) \circ \ad_\hbar(F_{1,\alpha}^{P'P_\lef(k')\Sigma\nu}\otimes
1^{\otimes k''})\cdots \ad_\hbar(F_{\nu,\alpha}^{P'P_\lef(k')\Sigma\nu}\otimes
1^{\otimes k''})\\
\circ ((\delta^{(R)}\otimes
\delta^{(P'')})\circ\delta)^{\Sigma,k'+1,\dots,n+1},
\end{multline*}
which is (\ref{delta:P:P0}) for $(P,P_1)$. In the same way, one
proves a similar identity relating $P_1$ and $P_2$. Let us now prove the
identity relating $P_2$ and $P_0$. We have
$\delta^{(P_2)}=(\delta\otimes \Id^{\otimes n-1})\circ \delta^{(P_2')}$ and
$\delta^{(P_0)}=(\delta\otimes \Id^{\otimes n-1})\circ \delta^{(P_0')}$, where
$P_2'$ and $P_0'$ are $n$-trees. We have
$$\Phi^{P_2,P_0}=(\Delta \otimes \Id^{\otimes n-2}) \circ
\Phi^{P_2',P_0'}$$
so we get
\begin{align*}
\delta^{(P_2)}&=\Ad(\Phi^{{P}_2,{P}_0})\circ \delta^{(P_0)}\\
&+ \big( \Ad(\Phi^{{P}_2,{P}_0})-\Ad((\Phi^{P_2',P_0'})^{1,3,\dots,n+1}) \big) \circ 
(\delta^{(P_0')})^{1,3,\dots,n+1}\\
&+ \big( \Ad(\Phi^{{P}_2,{P}_0})-\Ad((\Phi^{P_2',P_0'})^{2,3,\dots,n+1}) \big) \circ
(\delta^{(P_0')})^{2,3,\dots,n+1} \\
&+(\delta \otimes \Id^{\otimes n-1})\big(\sum_{k \leq n}\sum_{R\pa k\ptree}
\sum_{\begin{array}{l}
{\scriptstyle{\Sigma \subset\{1,\dots,n\},}}\cr
{\scriptstyle{\pecard(\Sigma)=k}}
\end{array}}\\
&\hskip1cm \sum_{\nu \geq 1} \sum_\alpha 
\ad_\hbar(F_{1,\alpha}^{P_2'P_0'\Sigma\nu})\cdots
\ad_\hbar(F_{\nu,\alpha}^{P_2'P_0'\Sigma\nu})\circ (\delta^{(R)})^\Sigma
\big).
\end{align*}
We have $\hbar \log \Phi^{P_2,P_0}\in {U'}^{\bo n+1}$ and $\hbar \log
\Phi^{P_2',P_0'}\in {U'}^{\bo n}$; this fact and the relations
\begin{multline*}
(\delta \otimes \Id^{\otimes n-1})(\ad_\hbar (x_1) \cdots \ad_\hbar (x_\nu)
\circ (\delta^{(R)})^\Sigma)
=\\
\big(\ad_\hbar(x_1^{12,\dots,n+1})\circ \cdots \circ
\ad_\hbar(x_\nu^{12,\dots,n+1})
- \ad_\hbar(x_1^{1,3,\dots,n+1})\circ \cdots \circ
\ad_\hbar(x_\nu^{1,3,\dots,n+1})\\
- \ad_\hbar(x_1^{2,3,\dots,n+1})\circ \cdots \circ
\ad_\hbar(x_\nu^{2,3,\dots,n+1})\big)\circ(\delta^{(R)})^{\Sigma +1}
\end{multline*}
if $1 \notin \Sigma$, and
\begin{align*}
&(\delta \otimes \Id^{\otimes n-1})(\ad_\hbar (x_1) \cdots \ad_\hbar (x_\nu)
\circ (\delta^{(R)})^\Sigma)=\\
&\hskip0.5cm\ad_\hbar(x_1^{12,\dots,n+1})\circ \cdots \circ
\ad_\hbar(x_\nu^{12,\dots,n+1})\circ ((\delta \otimes \Id^{\otimes n-1})\circ
\delta^{(R)})^{1,2,\Sigma'+1}\\
&\hskip0.5cm +\big( \ad_\hbar(x_1^{12,\dots,n+1})\circ \cdots \circ
\ad_\hbar(x_\nu^{12,\dots,n+1})
-\ad_\hbar(x_1^{1,3,\dots,n+1})\circ \cdots \circ
\ad_\hbar(x_\nu^{1,3,\dots,n+1})\big) \\
&\hskip10cm \circ(
\delta^{(R)})^{1,\Sigma'+1}\\
&\hskip0.5cm +\big( \ad_\hbar(x_1^{12,\dots,n+1})\circ \cdots \circ
\ad_\hbar(x_\nu^{12,\dots,n+1})
-\ad_\hbar(x_1^{2,3,\dots,n+1})\circ \cdots \circ
\ad_\hbar(x_\nu^{2,3,\dots,n+1})\big) \\
&\hskip10cm \circ(
\delta^{(R)})^{2,\Sigma'+1}.
\end{align*}
if $\Sigma=\Sigma' \cup \{1 \}$, where $1 \notin \Sigma'$, imply that
$\delta^{(P_2)}-\Ad(\Phi^{P_2,P_0})\circ \delta^{(P_0)}$ has
the desired form.
\smallskip
Let us now treat the case $|P_1|=1$. For this, we introduce the trees $P_3$ and
$P_4$, such that:
\begin{align*}
\delta^{(P_3)}&=(\Id^{\otimes n-1} \otimes \delta) \circ (\Id^{\otimes n-2}
\otimes \delta)\circ \cdots \circ \delta,\\
\delta^{(P_4)}&=(\Id^{\otimes n-1} \otimes \delta) \circ 
(\delta \otimes \Id^{\otimes n-2})\circ (\delta \otimes \Id^{\otimes n-3}
)\circ \cdots \circ (\delta \otimes \Id)\circ \delta.
\end{align*}
We then prove the relation for the pair $(P,P_3)$ in the same way as for
$(P_1,P_2)$ (only the right branch of the tree is changed);
the relation for $(P_3,P_4)$ in the same way as for $(P_2,P_3)$
(instead of composing a known relation by $\delta \otimes \Id^{\otimes n-1}$,
we compose it with $\Id^{\otimes n-1}\otimes \delta$); and using the identity
$$\delta^{(P_4)}=(\delta \otimes \Id^{\otimes n-1}) \circ
(\Id^{\otimes n-2} \otimes \delta) \circ (\delta \otimes \Id^{\otimes n-3})\circ
\cdots \circ \delta,$$
we prove the relation for $(P_4,P)$ in the same way as for $(P_2,P_3)$
(composing a known relation by $\delta \otimes \Id^{\otimes n-1}$).
\findem
\bigskip

\noindent {\bf - b - Properties of $\delta^{(P)}$}

\bigskip

\begin{proposition}
\label{prepa:flatness}
Let $n$ be an integer and $x \in U$.
\begin{description}
\item[1)] Assume that for any tree $R$, such that $|R|<n$, we have
$\delta^{(R)}(x)\in \hbar^{|R|}U^{\ho |R|}$. Then the conditions
\begin{equation}
\label{CP}
\delta^{(P)}(x)\in \hbar^n U^{\ho n}
\end{equation}
where $P$ is an $n$-tree, are all equivalent.
\item[2)] Assume that for any tree $R$, such that $|R|<n$, we have
$\delta^{(R)}(x) \in \hbar^{|R|+1}U^{\ho |R|}$. Then
the elements
$$\left({{1}\over{\hbar^n}}\delta^{(P)}(x)\mod \hbar\right) \in U(\g)^{\otimes
n},$$
where $P$ is an $n$-tree, are all equal and belong to $(\g^{\otimes
n})^{\fgam_n}=S^n(\g).$
\end{description}
\end{proposition}
\dem
Let us prove 1). We have $\delta^{(P)}=(\Id - \eta \circ \epsilon)^{\otimes
|P|}\circ \delta^{(P)}$, where $\eta : \KM[[\hbar]] \to U$ is the unit map of
$U$, so  
\begin{multline*}
\delta^{(P)}=\Ad(\Phi^{P,P_0})\circ \delta^{(P_0)}+ 
\sum_{k|k<n}\sum_{R \pa k\ptree}\sum_{\begin{array}{l}
{\scriptstyle{\Sigma \subset\{1,\dots,n\},}}\cr
{\scriptstyle{\pecard(\Sigma)=k}}
\end{array}}
\sum_{\nu \geq 1} \sum_{\alpha} \\
(\Id - \eta \circ \epsilon)^{\otimes n}\circ
\ad_\hbar
(F_{1,\alpha}^{PP_0R\Sigma \nu})\circ \cdots
\circ \ad_\hbar ( F_{\nu,\alpha}^{PP_0R\Sigma\nu})\circ(\delta^{(R)})^\Sigma.
\end{multline*}
Then 1) follows from:
\begin{lemma}  
Let $\Sigma$ be a subset of $\{1,\dots,n\}$ (we will write $|\Sigma|$ instead
of $\Card(\Sigma)$) and let $U_0$
be the kernel of the counit of $U$.
Let $x \in \hbar^{|\Sigma|}(U_0)^{\ho |\Sigma|}$ and $F_1,\dots,F_\nu$ be
elements of
$({U'})^{\bo n}$.
Then
$$(\Id - \eta \circ \epsilon)^{\otimes n}(\ad_\hbar(F_1)\cdots
\ad_\hbar(F_\nu)(x^\Sigma))\in \hbar^n(U_0)^{\ho n}.$$
\end{lemma}
\noindent{\sc Proof of Lemma.~}
Each element $F\in ({U'})^{\bo n}$ is uniquely expressed
as a sum $F=\sum_{\Sigma \in \PC(\{1,\dots,n\})} F_\Sigma$,
where $F_\Sigma$ belongs to the
image of
\begin{align*}
(U'_0)^{\bo |\Sigma|} &\to
(U')^{\bo n},\\ 
f& \mapsto f^\Sigma,
\end{align*}
$\PC(\{1,\dots,n\})$ is the set of subsets of $\{1,\dots,n\}$, 
and $U'_0$ is the kernel of the counit of $U'$.
Then
\begin{multline*}
(\Id -\eta \circ \epsilon)^{\otimes n}(\ad_\hbar (F_1) \cdots
\ad_\hbar(F_\nu)(x^\Sigma))\\
=\sum_{\Sigma_1,\dots,\Sigma_\nu \in \PC(\{1,\dots,n\})}(\Id -\eta \circ
\epsilon)^{\otimes n}
\big(\ad_\hbar((F_1)_{\Sigma_1})\cdots
\ad_\hbar((F_\nu)_{\Sigma_\nu})(x^\Sigma)\big).
\end{multline*}
The summands corresponding to $(\Sigma_1,\dots,\Sigma_\nu)$ such that
$\Sigma_1 \cup \cdots \Sigma_\nu \cup \Sigma \not=\{1,\dots,n\}$ are all zero.
Moreover, each 
$(F_\alpha)_{\Sigma_\alpha}$ can be expressed as $(f_\alpha)^{\Sigma_\alpha}$,
where
$f_\alpha \in \hbar^{|\Sigma_\alpha|}(U_0)^{\ho |\Sigma_\alpha|}$.
The lemma then follows from the statement:
\begin{statement}
If $\Sigma,\Sigma'\subset \{1,\dots,n\}$, $x \in \hbar^{|\Sigma|}(U_0)^{\ho
|\Sigma|}$,
$y \in \hbar^{|\Sigma'|}(U_0)^{\ho |\Sigma'|}$, then
${{1}\over{\hbar}}[x,y]$ can be expressed as $z^{\Sigma \cup \Sigma'}$, where
$z\in \hbar^{|\Sigma \cup \Sigma'|} (U_0)^{\ho |\Sigma \cup \Sigma'|}.$
\end{statement}
\dem
If $\Sigma \cap \Sigma' = \emptyset$, then $[x,y]=0$,
so the statement holds. If
$\Sigma \cap \Sigma'\not=\emptyset$, then the
$\hbar$-adic valuation
of ${{1}\over{\hbar}}[x,y]$ is
$\geq -1+|\Sigma|+|\Sigma'|\geq |\Sigma|+|\Sigma'|-|\Sigma \cap \Sigma'|=|\Sigma
\cup \Sigma'|$.
\findem

\smallskip
Let us now prove property 2). The above arguments immediately imply that the
$({{1}\over{\hbar^n}}\delta^{(P)}(x) \mod \hbar)$, $|P|=n$, are all
equal. This defines an element $S_n(x)\in U(\g)^{\otimes n}.$ If $|P|=n$, we
have 
$(\Id^{\otimes k}\otimes \delta \otimes \Id^{\otimes n-k-1})\circ
\delta^{(P)}(x)\in
\hbar^{n+1}U^{\ho n+1}$, so if $\delta_0$~:
$U(\g) \to U(\g) \otimes U(\g)$ is defined by
$\delta_0(x)=\Delta_0(x)-x\otimes 1 -1\otimes x+ \epsilon(x) 1 \otimes 1$,
$\Delta_0$ being the
coproduct of $U(\g)$,
then
$(\Id^{\otimes k}\otimes \delta_0 \otimes \Id^{\otimes n-k-1})(S_n(x))=0$,
so \begin{equation}
\label{s:g:n}
S_n \in \g^{\otimes n}.
\end{equation}
Let us denote by
$\sigma_{i,i+1}$ the permutation of the factors $i$ and $i+1$ in a tensor
power. For $i=1,\dots,n-1$, let us
compute
$(\sigma_{i,i+1}-\Id)(S_n(x))$. Let $P'$ be
a $(n-1)$-tree and let $P$ be the
$n$-tree such that 
$\delta^{(P)}=(\Id^{\otimes i-1} \otimes \delta \otimes \Id^{\otimes
n-i-1})\circ \delta^{(P')}$.
Then
$$(\sigma_{i,i+1}-\Id)(S_n)=\left[{{1}\over{\hbar}}(\Id^{\otimes i-1}
\otimes
(\delta^{2,1}-\delta)\otimes \Id^{\otimes n-i-1})\circ \delta^{(P')}(x)\mod
\hbar\right].$$
By assumption, $\delta^{(P')}(x) \in \hbar^n U^{\ho n-1}$; moreover,
$\delta^{2,1}-\delta=\Delta^{2,1} - \Delta$, so
$(\delta^{2,1}-\delta)(U) \subset \hbar (U\ho U)$; therefore
$$(\Id^{\otimes i-1}\otimes (\delta^{2,1}-\delta)\otimes \Id^{\otimes
n-i-1})\circ \delta^{(P')}(x) \in \hbar^{n+1} U^{\ho n};$$
it follows that $(\sigma_{i,i+1}-\Id)(S_n(x))=0$, therefore $S_n(x)$ is a
symmetric tensor of $U(\g)^{\otimes n}$. Together with
(\ref{s:g:n}), this gives $S_n(x)\in (\g^{\otimes n})^{\fgam_n}.$ This ends
the proof of Proposition \ref{prepa:flatness}.
\findem

\bigskip

\noindent {\bf - c - Flatness of $U'$ (proof of Theorem \ref{theo:2})}

\bigskip

Let us set
$${U''}^{(n)}=\{x \in U'|\delta^{(P)}(x) \in \hbar^{|P|+1}U^{\ho |P|}\hbox{ if
}|P|\leq n\}.$$
Then by Proposition \ref{multiplication}, we have a decreasing algebra
filtration
\begin{equation}
\label{decr:fil}
U'={U''}^{(0)}\supset {U''}^{(1)} \supset {U''}^{(2)} \supset \cdots
\supset \hbar U'.
\end{equation}
Each ${U''}^{(n)}$ is divisible in $U'$, i.e., ${U''}^{(n)} \cap \hbar U' =
\hbar {U''}^{(n)}$.
We also
have ${U''}^{(n)}\supset {U'}^{(n)} + \hbar U'$ (we will see later that this is an
equality). We derive from (\ref{decr:fil}) a decreasing filtration
$$\OC={\OC''}^{(0)} \supset {\OC''}^{(1)} \supset {\OC''}^{(2)} \supset
\cdots,$$
where $\OC=U'/\hbar U'$ and ${\OC''}^{(n)}={U''}^{(n)}/\hbar {U''}^{(n)}$.
We have clearly
$$\CAP_{n \geq 0}{\OC''}^{(n)}=\{0\};$$
the fact that $\OC$ is complete for this filtration will follow from its
identification
with the filtration $\OC \supset {\OC'}^{(1)}\supset \cdots $ 
(see Proposition \ref{completeness}), where ${\OC'}^{(i)}={U'}^{(i)}/\hbar U
\cap {U'}^{(i)}$ and ${U'}^{(i)}$ is
defined in (\ref{filtr:1}).
We first prove:
\begin{proposition}
Set $\widehat{\gr}''(\OC)=\hopl_{n \geq 0} {\OC''}^{(n)}/{\OC''}^{(n+1)}$.
Then there is a unique
linear map
$\lambda_n$~: ${\gr}_n''(\OC) \to S^n(\g)$, taking the class of $x$ to the
common value of
all
${{1}\over{n!}}({{1}\over{\hbar^n}} \delta^{(P)}(x)\mod \hbar)$, where $P$ is a
$n$-tree.
The resulting map
$\lambda~:~\widehat{\gr}''(\OC)\to{\widehat{S}}^\cdot (\g)$ is an isomorphism of
graded complete algebras.
\end{proposition}
\dem In Proposition \ref{prepa:flatness}, we constructed a map 
${U''}^{(n)}\to S^n(\g)$, by $x \mapsto $ common value of 
${{1}\over{n!}}({{1}\over{\hbar^n}} \delta^{(P)}(x)\mod \hbar)$ for all
$n$-trees $P$.
The subspace ${U''}^{(n+1)} \subset {U''}^{(n)}$ is clearly contained in the
kernel of this map, so we obtain a map
$$\lambda_n~:~{U''}^{(n)}/{U''}^{(n+1)}={\OC''}^{(n)}/{\OC''}^{(n+1)}\to
S^n(\g).$$
Let us prove that $\lambda=\hopl_{n \geq 1}\lambda_n$ is a morphism of algebras.
If $x \in {U''}^{(n)}$ and
$y \in {U''}^{(m)}$, Proposition \ref{multiplication} implies that if
$R$ is any $(n+m)$-tree, we have
$$\delta^{(P)}(xy)=\sum_{\begin{array}{l}
{\scriptstyle{\Sigma',\Sigma''\subset\{1,\dots,n+m\}| }}\cr
{\scriptstyle{\Sigma' \cup \Sigma''=\{1,\dots,n+m\}}}
\end{array}}
\delta^{(R_{\Sigma'})}(x)^{\Sigma'}\delta^{(R_{\Sigma''})}(y)^{\Sigma''}.$$
The $\hbar$-adic valuation of the term corresponding to $(\Sigma',\Sigma'')$
is $\geq |\Sigma'|+|\Sigma''|$ if
$|\Sigma'|\geq n$ and
$|\Sigma''|\geq m$,
and $\geq \!|\Sigma'|\! + \!|\Sigma''|\!+\!1$ otherwise,
so the only contributions to
$(\!{{1}\over{\hbar^{n+m}}}\delta^{(R)}(xy) \mod \hbar)$ are
those
of the pairs $(\Sigma',\Sigma'')$ such that $\Sigma' \cap \Sigma'' =\emptyset$.
Then:
\begin{align*}
&({{1}\over{\hbar^{n+m}}}\delta^{(R)}(xy) \mod \hbar)\\
&\hskip1.5cm=
\sum_{\begin{array}{l}
{\scriptstyle{\Sigma',\Sigma''\subset\{1,\dots,n+m\}| }}\cr
{\scriptstyle{|\Sigma'|=n,|\Sigma''|=m,}}\cr
{\scriptstyle{\Sigma' \cap \Sigma''=\emptyset}}
\end{array}}
({{1}\over{\hbar^n}}\delta^{(R_{\Sigma'})}(x) \mod \hbar)
({{1}\over{\hbar^m}}\delta^{(R_{\Sigma''})}(y) \mod \hbar)\\
&\hskip1.5cm=\sum_{\begin{array}{l}
{\scriptstyle{\Sigma',\Sigma''\subset\{1,\dots,n+m\}| }}\cr
{\scriptstyle{|\Sigma'|=n,|\Sigma''|=m,}}\cr
{\scriptstyle{\Sigma' \cap \Sigma''=\emptyset}}
\end{array}}
(n!\lambda_n(x)^{\Sigma'})(m!\lambda_m(y)^{\Sigma''})\cr
&\hskip1.5cm=(n+m)!\lambda_n(x)\lambda_m(y),
\end{align*}
because the map 
\begin{align*}
S^\cdot(\g) &\to (T(\g),\hbox{shuffle product}),\\
x_1\cdots x_n &\mapsto 
\sum_{\sigma \in \fgam_n} x_{\sigma(1)}\otimes \cdots \otimes
x_{\sigma(n)}
\end{align*}
is an algebra morphism.
Therefore
$\lambda_{n+m}(xy)=\lambda_n(x)\lambda_m(y)$.
Let us prove that $\lambda_n$ is injective. If $x \in {U''}^{(n)}$ is such that
$({{1}\over{\hbar^n}}\delta^{(P)}(x)\mod \hbar)=0$ for any $n$-tree $P$, then
$x \in {U''}^{(n+1)}$,
so its class in ${\OC''}^{(n)}/{\OC''}^{(n+1)}={U''}^{(n)}/{U''}^{(n+1)}$ is
zero.
So each $\lambda_n$ is injective, so $\lambda$ is injective.

\smallskip

\noindent To prove that $\lambda$ is surjective, it suffices to prove that
$\lambda_1$ is surjective. Let us fix $x \in \g$. We will construct a 
sequence $x_n \in U$, $n\geq 0$ such that $\epsilon(x_n)=0$,
$({{1}\over{\hbar}} x_n ~\mod~ \hbar)=x$,
$x_{n+1} \in x_n + \hbar^{n+1} U$ for any $n \geq 1$, and if $P$ is any
tree such that $|P| \leq n$, $\delta^{(P)}(x_n) \in
\hbar^{|P|} U^{\ho|P|}$ (this last
condition
implies that
$\delta^{(Q)}(x_n)\in \hbar^n U^{\ho |Q|}$ for
$|Q| \geq n$).
Then the limit
$\widetilde{x}=\Lim_{n \to \infty}(x_n)$ exists, belongs to $U'$, satisfies
$\epsilon(\widetilde{x})=0$
and $({{1}\over{\hbar}} \delta_1(\widetilde{x}) \mod \hbar)=x$, so
its class in ${U''}^{(1)}/{U''}^{(2)}$ is a preimage of $x$.

\smallskip

\noindent Let us now construct the sequence $(x_n)_{n \geq 0}$. We fix a linear map
$\g \to \{y \in U|\epsilon(y)=0\},~y\mapsto \bar{y}$, such that
for any $y\in \g$, $(\bar{y} \mod \hbar)=y$.
We set $x_1=\hbar \bar{x}$. Let us
construct
$x_{n+1}$ knowing $x_n$.
By Proposition \ref{prepa:flatness}, if $Q$ is any $(n+1)$-tree,
$\delta^{(Q)}(x_n) \in \hbar^n U^{\ho n+1}$, and
$({{1}\over{\hbar^n}}\delta^{(Q)}(x_n) \mod \hbar )$ is an element
of $S^{n+1}(\g)$, independent of $Q$. Let us write this element
as
$$\sum_{\sigma \in \fgam_{n+1}} \sum_\alpha y_{\sigma(1)}^\alpha \cdots 
y_{\sigma(n+1)}^\alpha,\hbox{ where }\sum_\alpha y_{1}^\alpha \otimes \cdots \otimes
y_{n+1}^\alpha \in \g^{\otimes n+1}.$$
Then we set
$$\hskip2.7cm x_{n+1}=x_n - {{\hbar^{n+1}}\over{(n+1)!}}\sum_{\sigma \in
\fgam_{n+1}}\bar{y}_{\sigma(1)}^\alpha\cdots
\bar{y}_{\sigma(n+1)}^\alpha.\hskip2.7cm\findem$$
We now prove:
\begin{proposition}
\label{completeness}
$\hbox{~}$

\begin{description}
\item[1)] For any $n \geq 0$, ${U''}^{(n)}={U'}^{(n)}+\hbar U';$
\item[2)] The filtrations ${\OC}={\OC'}^{(0)} \supset 
{\OC'}^{(1)} \supset \cdots$ and ${\OC}={\OC''}^{(0)} \supset {\OC''}^{(1)}
\supset \cdots$
coincide, and $\OC$ is complete and separated for this filtration.
\end{description}
\end{proposition}
\dem Let us prove 1). We have to show that ${U''}^{(n)} \subset
{U'}^{(n)}+\hbar U'$.
Let $x \in {U''}^{(n)}$. We have $\delta^{(P)}(x)\in \hbar^{|P|+1}U^{\ho |P|}$
for
$|P| \leq n-1$, and for $P$ an $n$-tree, $({{1}\over{\hbar^n}}\delta^{(P)}(x)
\mod \hbar) \in
S^n(\g)$ and is independent on $P$. Write this element of $S^n(\g)$ as
$\sum_{\sigma \in \fgam_n} \sum_\alpha y_{\sigma(1)}^\alpha \otimes \cdots 
y_{\sigma(n)}^\alpha$ and set
$f_n={{1}\over{n!}}\sum_{\sigma \in \fgam_n} \sum_\alpha 
\bar{y}_{\sigma(1)}^\alpha \cdots 
\bar{y}_{\sigma(n)}^\alpha$. Then each $\bar{y}_{i}^\alpha$ belongs to
$U'\cap \hbar U$, so
$f_n \in U' \cap \hbar^n U = {U'}^{(n)}$.
Moreover, $x - f_n$ belongs to
${U''}^{(n+1)}$. Iterating this procedure, we construct elements
$f_{n+1},f_{n+2},\dots,$ where each $f_k$ belongs to
${U'}^{(k)}$.
The series
$\sum\limits_{k \geq n} f_k$ converges in $U'$; denote by $f$ its sum, then
$x-f$ belongs to
$\CAP_{k \geq n} {U''}^{(k)}=\hbar U'$. So
${U''}^{(n)} \subset {U'}^{(n)} + \hbar U'$.
The inverse inclusion is obvious.
This proves 1). Then 1) immediately implies that for any $n$, ${\OC'}^{(n)}={\OC''}^{(n)}$.
We already know $\OC$ is complete and separated for $\OC={\OC'}^{(0)}\supset
{\OC'}^{(1)}\supset \cdots$, which proves 2).
$\hbox{~}$\findem

\medskip

\noindent{\sc End of proof of Theorem \ref{theo:2}.~}
$\OC$ is a complete local ring, and we have
a ring isomorphism $\widehat{\gr}(\OC) \to \widehat{S}^\cdot(\g)$.
Then any lift
$\g \to {\OC'}^{(1)}$ of ${\OC'}^{(1)} \to {\OC'}^{(1)}/{\OC'}^{(2)}=\g$ yields a continuous
ring morphism $\mu$~: $\widehat{S}^\cdot(\g) \to \OC$.
The associated graded of $\mu$ is the identity, so $\mu$ is an isomorphism.
So $\OC$ is noncanonically isomorphic to $\widehat{S}^\cdot(\g)$.
\findem
\begin{remark}
\label{noncanonicity}
When $U$ is Hopf and $\g$ is finite-dimensional, $U'/\hbar U'$ identifies
canonically with $\OC_{G^*}=(U(\g^*))^*$, where $\g^*$ is the dual Lie
bialgebra
of $\g$ (see \cite{Dr1}, \cite{Ga}).
The natural projection ${T}(\g^*) \to U(\g^*)$ and
the identification ${T}(\g^*)^*=\widehat{T}(\g)$ (where $\widehat{T}(\g)$ means
the degree completion) induce an injection $U'/\hbar U'=\OC_{G^*}=(U(\g^*))^*
\hookrightarrow \widehat{T}(\g)$.
The map
$U'/\hbar U'
\hookrightarrow \widehat{T}(\g)$ can be interpreted simply as follows.
For any
$x \in U'$, we have $({{1}\over{\hbar^n}}\delta_n(x) \mod \hbar) \in \g^{\otimes n}$.
Then $U'/\hbar U'
\hookrightarrow \widehat{T}(\g)$ takes the class of $x \in U'$ to the sequence 
$({{1}\over{\hbar^n}}\delta_n(x) \mod \hbar)_{n \geq 0}$.

In the quasi-Hopf case, we have no canonical embedding
$U'/\hbar U'
\hookrightarrow \widehat{T}(\g)$ because the various $({{1}\over{\hbar^n}}\delta^{(P)}(x) \mod
\hbar)$ do not necessarily coincide for all the $n$-trees $P$.
This is related to the
fact that one cannot expect a Hopf pairing
$U(\g^*)\otimes (U'/\hbar U')\to \KM$ since $\g^*$ is no longer a Lie algebra,
so $U(\g^*)$ does not make sense.

In the other hand, Theorem \ref{theo:2} can be interpreted as follows: in the
Hopf case, the exponential induces an isomorphism of formal schemes
$\g^* \to G^*$, so $U'/\hbar U'$ identifies noncanonically with
$\OC_{\g^*}=\widehat{S}^\cdot(\g)$.
In the quasi-Hopf case, although there is no formal group $G^*$, we still have
an isomorphism $U'/\hbar U' \To^{\sim} \widehat{S}^\cdot(\g)$.
\findem 
\end{remark}

\vskip20pt

\centerline {\bf \S\; 4 \ Twists}

\stepcounter{section}\label{section:4}

\vskip20pt

\noindent {\bf - a - Admissible twists}

\bigskip

If $(U,m,\Delta,\Phi) $ is an arbitrary QHQUE algebra, we will call
a twist $F \in (U^{\ho 2})^\times$ {\it admissible} if
$\hbar \log(F)\in (U')^{\bo 2}$.
\begin{proposition}
\label{adm:twist}
Let $(U,m,\Delta,\Phi)$ be an admissible quasi-Hopf algebra and
$F$ an admissible twist. Then the twisted
quasi-Hopf algebra $(U,m,\!{}^F\!\Delta,\!{}^F\!\Phi)$ is admissible.
\end{proposition}
\dem
Set $f = \hbar \log(F)$. Then we have
$$\hbar \log({}^F\Phi)=f^{1,2} \star f^{12,3} \star
(\hbar \log(\Phi))\star(-f^{1,23})\star(-f^{2,3}),$$
where $a \star b = a+b + {{1}\over{\hbar}} [a,b] + \cdots$ (the CBH series for
${U'}^{\bo 3}$ equipped with the bracket ${{1}\over{\hbar}}[-,-]$).
Since ${U'}^{\bo 3}$ is stable under $\star$, we have $\hbar \log({}^F\Phi)
\in {U'}^{\bo 3}$. So
$(U,m,{}^F\Delta,{}^F\Phi)$ is admissible.
\findem

\medskip

Let us now prove 
\begin{proposition}
\label{twist:egal}
Under the hypothesis of Proposition \ref{adm:twist},
the QFS algebra ${U}_F'$ corres-ponding to 
$(U,m,\!{}^F\!\!\Delta,\!{}^F\!\Phi)$ coincides with the QFS algebra $U'$
corresponding to $(U,m,\Delta,\Phi)$.
\end{proposition}
\noindent We will first prove the following lemma:
\begin{lemma}
\label{variation:F}
Let $P$ be an $n$-tree. Then
\begin{align}
\label{eq:*}
\begin{split}
\delta_F^{(P)}&=\delta^{(P)} + \sum_{k \leq n}
\sum_{R\pa k\ptree}
\sum_{\begin{array}{l}
{\scriptstyle{\Sigma\subset\{1,\dots,n\}| }}\cr
{\scriptstyle{\pecard(\Sigma)=k}}\cr
\end{array}}\\
&\hskip4cm \sum_{\nu \geq 1} \sum_\alpha \ad_\hbar(f_{1,\alpha}^{\Sigma,P})\circ \cdots\circ
\ad_\hbar(f_{\nu,\alpha}^{\Sigma,P})\circ(\delta^{(R)})^\Sigma,
\end{split}
\end{align}
where for each $\nu$, $\sum_\alpha f_{1,\alpha}^{\Sigma,P} \otimes \cdots \otimes
f_{\nu,\alpha}^{\Sigma,P} \in ({U'}^{\bo n})^{\bo \nu}$.
\end{lemma}
\begin{remark}
One can prove that in the right hand side of (\ref{eq:*}), the contribution 
of all terms with $k=n$ is
$(\Ad(F^{(P)})-\Id)\circ \delta^{(P)}$ where $F^{(P)}$ is
the product of
$F^{I,J}$ ($I,J$ subsets of $\{1,\dots,n\}$, such that $\max(I) < \min(J)$) and
their inverses such that
$$\Delta_F^{(P)}=\Ad(F^{(P)}) \circ \Delta^{(P)}.$$
\end{remark}
\noindent{\sc Proof of the lemma.~}
equation (\ref{eq:*}) may be proved by induction on $|P|$. Let us prove it for
the unique tree $P$ such that $|P|=2$:
$$\delta_F^{(2)}=\delta^{(2)}+ \sum_{\nu \geq 1}{{1}\over{\nu !}} \ad_\hbar(f)^\nu
(\delta^{(2)}(x)+\delta^{(1)}(x)^1 + \delta^{(1)}(x)^2),$$
where $(1)$ and $(2)$ are the $1$- and $2$-trees.
Assume that (\ref{eq:*}) is proved when $|P|=n$. Let
$P'$ be an $(n+1)$-tree. Then for some $i\in \{1,\dots,n\}$, we have
$$\delta_F^{(P)}=(\Id^{\otimes i-1} \otimes \delta_F^{(2)} \otimes
\Id^{\otimes n-i})\circ \delta_F^{(P')},$$
where $|P'|=n$. Then:
\begin{align*}
\delta_F^{(P)}=&(\Id^{\otimes i-1}\otimes \Delta_F \otimes \Id^{\otimes
n-i})\circ
\delta_F^{(P')} -(\delta_F^{(P')})^{1,\dots,\widehat{i},\dots, n+1}-
(\delta_F^{(P')})^{1,\dots,\widehat{i+1},\dots, n+1}\\
=&(\Id^{\otimes i-1}\otimes \Delta_F \otimes \Id^{\otimes
n-i})\circ
\big(\delta^{(P')} + \sum_{k \leq n}\sum_{R\pa k\ptree}
\sum_{\begin{array}{l}
{\scriptstyle{\Sigma\subset\{1,\dots,n\}| }}\cr
{\scriptstyle{\pecard(\Sigma)=k}}\cr
\end{array}}\\
&\hskip4.4cm 
\sum_{\nu \geq 1} \sum_\alpha \ad_\hbar(f_{1,\alpha}^{\Sigma,P'})\circ 
\cdots\circ
\ad_\hbar(f_{\nu,\alpha}^{\Sigma,P'})\circ(\delta^{(R)})^\Sigma\big)\\
&-\big(\cdots \big)^{1\dots
\widehat{i},\dots, n+1}-\big(\cdots\big)^{1,\dots, 
\widehat{i+1} ,\dots, n+1}\\
=&\Ad(F^{i,i+1})\circ\big(\delta^{(P)}+ (\delta^{(P')})^{1,\dots, \widehat{i},\dots,
n+1}+(\delta^{(P')})^{1,\dots, \widehat{i+1},\dots, n+1}\\
&+ \sum_{k \leq n }\sum_{R\pa k\ptree}
\sum_{\begin{array}{l}
{\scriptstyle{\Sigma\subset\{1,\dots,n\}| }}\cr
{\scriptstyle{\pecard(\Sigma)=k}}\cr
\end{array}}\sum_{\nu \geq 1} \sum_\alpha 
\ad_\hbar ((f_{1,\alpha}^{\Sigma,P'})^{1,\dots,\{i,i+1\},\dots,n+1})\circ\\
&\hskip2cm \circ
\ad_\hbar((f_{\nu,\alpha}^{\Sigma,P'})^{1,\dots,\{i,i+1\},\dots,n+1})\circ
(1^{\otimes i-1}\otimes \Delta \otimes 1^{\otimes n-i})
\circ(\delta^{(R)})^\Sigma{\big)}\\
&-\big(\cdots\big)^{1,\dots, \widehat{i} ,\dots, n+1}- 
\big(\cdots\big)^{1,\dots, \widehat{i+1} ,\dots, n+1};
\end{align*}
this has the desired form because:
\begin{multline*}
(\Ad(F^{i,i+1})-1)\circ \big(\delta^{(P)}+(\delta^{(P')})^{1,\dots, \widehat{i}
,\dots, n+1} + (\delta^{(P')})^{1 ,\dots, \widehat{i+1} ,\dots, n+1}\big)\\
=\sum_{\nu \geq 1} {{1}\over{\nu !}} \ad_\hbar(f^{i,i+1})^\nu \big(\delta^{(P)}+
(\delta^{(P')})^{1,\dots, \widehat{i},\dots, n+1}+ (\delta^{(P')})^{1,\dots,
\widehat{i+1},\dots, n+1}\big).
\end{multline*}
This proves (\ref{eq:*}).
\findem

\medskip

\noindent{\sc End of proof of Proposition \ref{twist:egal}.~}
One repeats the proof of Proposition \ref{prepa:flatness} to
prove that if $x \in U'$, then we have $\delta^{(P)}(x) \in
\hbar^{|P|}U^{\ho |P|}$ for any tree $P$.
So $U' \subset {U}_F'$.
Since $(U,m,\Delta,\Phi)$ is the twist by $F^{-1}$ of
$(U,m,{}^F\!\Delta,{}^F\!\Phi)$, and
$\hbar \log(F^{-1})=-\hbar \log(F) \in (U')^{\bo 2} \subset (U_F')^{\bo 2}$,
$F^{-1}$ is admissible for
$(U,m,{}^F\!\Delta,{}^F\!\Phi)$, so we have also ${U}_F' \subset U'$, so
${U}_F'=U'$.
\findem

\bigskip

\noindent {\bf - b - Twisting any algebra into
an admissible algebra}

\bigskip

\begin{proposition}
\label{adm:twist2}
Let $(U,m,\Delta,\Phi)$ be a quasi-Hopf algebra. There exists a twist $F_0$
such that the twisted
quasi-Hopf algebra $(U,m,{}^{F_0}\Delta,{}^{F_0}\Phi)$ is admissible.
\end{proposition}
\dem
We construct $F_0$ as a convergent infinite product
$F_0= \cdots F_n \cdots F_2$, where $F_n \in 1 + \hbar^{n-1}U^{\ho 2}$, and the
$F_n$ have the
following property: if ${\bar{F}}_n=F_nF_{n-1} \cdots F_2$, if
$\Phi_n={}^{{\bar{F}}_n}\Phi$, and
$\delta_n^{(P)}$~: $U \to U^{\ho |P|}$ is the map corresponding to a tree $P$ and
to
$\Delta_n = \Ad({\bar{F}}_n)\circ \Delta$, then we have
$$(\delta_n^{(P)} \otimes
\delta_n^{(Q)} \otimes \delta_n^{(R)})
(\hbar \log (\Phi_n)) \in
\hbar^{|P| + |Q| + |R|}
U^{\ho |P| + |Q| + |R|}$$
for any trees $P,Q,R$ such that $|P| + |Q| + |R|\leq n$.

\medskip

Assume that we have constructed $F_1,\dots,F_n$, and let us construct $F_{n+1}$.
The argument
of Proposition \ref{prepa:flatness} shows that for any integers $(n_1,n_2,n_3)$
such
that $n_1+n_2+n_3=n+1$, and
any trees $P,Q,R$ such that $|P|=n_1$, $|Q|=n_2$, $|R|=n_3$,
$$\big({{1}\over{\hbar^n}}(\delta_n^{(P)}\otimes
\delta_n^{(Q)}\otimes \delta_n^{(R)})(\hbar \log(\Phi_n))\mod \hbar\big) \in
S^{n_1}(\g)\otimes S^{n_2}(\g) \otimes S^{n_3}(\g),$$
and is independent of the trees $P$, $Q$, $R$.
The direct sum of these elements is an element $\bar{\varphi}_n$ of $ S^\cdot(\g)^{\otimes 3}$,
homogeneous of degree $n+1$. 
Since $\Phi_n$ satisfies the pentagon equation
$$(\Id \!\otimes \!\Id \!\otimes \!\Delta_n)(\Phi_n\!)^{-1}\!(1 \otimes \Phi_n)(\Id \otimes
\Delta_n \otimes \Id)(\Phi_n\!)(\Phi_n \otimes 1)(\Delta_n \otimes \Id \otimes
\Id)(\Phi_n)^{-1}
\!\!\!=1,$$
$\varphi_n^\hbar:=\hbar \log(\Phi_n)$ satisfies the equation
\begin{align}
\label{pent:phi:h}
\begin{split}
&\big(-(\Id \otimes \Id \otimes \Delta_n)(\varphi_n^\hbar)\big)\star
(1 \otimes \varphi_n^\hbar)\star \big((\Id \otimes \Delta_n \otimes \Id)
(\varphi_n^\hbar)\big)\star \\
&\hskip5cm (\varphi_n^\hbar \otimes 1) \star \big(-(\Delta_n \otimes \Id \otimes
\Id)(\varphi_n^\hbar)\big)=0,
\end{split}
\end{align}
where we set
$$a \star b = a+b + {{1}\over{2}} [a,b]_\hbar + \cdots$$
(the CBH series for the Lie bracket $[-,-]_\hbar$). Let
$(n_1,n_2,n_3,n_4)$ be integers such that $n_1 + \cdots + n_4=n+1$. Let
$P,Q,R,S$ be trees such that $|P|=n_1,\dots,|S|=n_4$. Let
us apply $\delta_n^{(P)} \otimes \cdots \otimes \delta_n^{(S)}$ to
(\ref{pent:phi:h}).
The left hand side of (\ref{pent:phi:h}) is equal to
\begin{equation*}
(-\Delta_n \otimes \Id \otimes \Id + \Id \otimes \Delta_n \otimes \Id
-\Id \otimes \Id \otimes \Delta_n)(\varphi_n^\hbar)
+ (1 \otimes \varphi_n^\hbar)-(\varphi_n^\hbar \otimes 1) + \hbox{ brackets}.
\end{equation*}
Now
\begin{equation*}
(\delta_n^{(P)} \otimes \delta_n^{(Q)} \otimes \delta_n^{(R)} \otimes
\delta_n^{(S)})(\Delta_n \otimes \Id \otimes \Id)(\varphi_n^\hbar)
=(\delta_n^{(P \cup Q)}\otimes \delta_n^{(R)} \otimes \delta_n^{(S)}).
\end{equation*}
where $P \cup Q$ is the tree with left descendant $P$ and right descendant $Q$.
Therefore
\begin{equation*}
\big({{1}\over{\hbar^n}}(\delta_n^{(P)} \otimes
\delta_n^{(Q)} \otimes \delta_n^{(R)} \otimes \delta_n^{(S)})(\Delta_n \otimes
\Id \otimes \Id)(\varphi_n^\hbar)\mod \hbar \big) 
= (\Delta_0 \otimes \Id \otimes \Id)(\bar{\varphi}_n)_{n_1,n_2,n_3,n_4}
\end{equation*}
where the index $(n_1,\dots,n_4)$ means the component in
$\Otimes_{i=1}^4 S^{n_i}(\g).$ On
the other hand, if $a_1$ and $a_2 \in U^{\ho 4}$ are such that
$$(\delta_n^{(P)} \otimes \cdots \otimes \delta_n^{(S)})(a_i) \in
\hbar^{\inf(|P|+\cdots + |S|,n)}U^{\ho 4}$$
for any trees $(P,\dots,S)$, then if $(P,\dots,S)$ are such that $|P|+ \cdots +
|S|=n$,
we have
$$(\delta_n^{(P)} \otimes \cdots \otimes
\delta_n^{(S)})({{1}\over{\hbar}}[a_1,a_2])\in \hbar^{n+1} U^{\ho n};$$
one proves this in the same way as the commutativity of $U'/\hbar U'$ (see
Theorem \ref{theo:1}). Then
${{1}\over{\hbar^n}}(\delta_n^{(P)}\otimes \cdots \otimes \delta_n^{(S)})
(\ref{pent:phi:h})|_{\hbar=0}$ yields 
$\dif(\bar{\varphi}_n)=0$, where $\dif$~:
$S^\cdot(\g)^{\otimes 2} \to S^\cdot(\g)^{\otimes 3}$ is the co-Hochschild
cohomology differential. This relation implies that
$$\bar{\varphi}_n=\dif(\bar{f}_n)+\lambda_n,$$
where $\bar{f}_n \in S^\cdot(\g)^{\otimes 2}$ and
$\lambda_n \in \Lambda^3(\g)$. Moreover, $f_n$ and $\lambda_n$ both have degree
$n+1$. This implies that $\lambda_n=0$.
Let $f_n \in (U(\g)^{\otimes 2})_{\leq n+1}$ be a preimage of $\bar{f}_n$ by
the projection 
$$(U(\g)^{\otimes 2})_{\leq n+1}\to (U(\g)^{\otimes 2})_{\leq n+1}/(U(\g)^{\otimes 2})_{\leq n}
= (S^\cdot(\g)^{\otimes 2})_{n+1}$$ 
(where the indices $n$ and $\leq n$ mean
``homogeneous part of degree $n$'' and ``part of degree $\leq n$'').
Let
$f_n^\hbar \in U^{\ho 2}$ be
a preimage of $f_n$ by the projection
$U^{\ho 2} \to U^{\ho 2}/\hbar U^{\ho 2}=U(\g)^{\otimes 2}$. Set
$F_{n+1}=\exp(\hbar^{n-1}f_n)$. We may assume that
$\hbar^n f_n \in (U(\bar{F}_n)')^{\bo 2}$, where $U(\bar{F}_n)'=
\{x \in U|\delta_n^{(P)}(x) \in
\hbar^{\inf(n,|P|)}U^{\ho |P|}\}$.
Then $\Phi_{n+1}={}^{F_{n+1}}\Phi_n$. If $P,Q,R$ are such that $|P|+|Q|+|R|=n+1$,
then 
$$(\delta_n^{(P)}\otimes \delta_n^{(Q)} \otimes \delta_n^{(R)})
(\hbar \log(\Phi_{n+1})) \in
\hbar^{n+1} U^{\ho n+1}.$$
Then according to Lemma \ref{variation:F},
$$(\delta_{n+1}^{(P)} \otimes \delta_{n+1}^{(Q)} \otimes \delta_{n+1}^{(R)} -
\delta_n^{(P)} \otimes \delta_n^{(Q)} \otimes \delta_n^{(R)})(\hbar
\log(\Phi_{n+1}))$$
has $\hbar$-adic valuation
$> |P|+|Q|+|R|$ when $|P|+|Q|+|R| \leq n+1$.
So
$(\delta_{n+1}^{(P)} \otimes \delta_{n+1}^{(Q)} \otimes \delta_{n+1}^{(R)})
(\hbar \log(\Phi_{n+1}))\in \hbar^{|P|+|Q|+|R|} U^{\ho |P|+|Q|+|R|}$
whenever $|P|+|Q|+|R| \leq n+1$.
\findem

\vskip20pt

\centerline {\bf \S\; 5 \ Proof of Proposition \ref{prop:1.4}}

\stepcounter{section}\label{section:5}

\vskip20pt

We work by induction on $n$.
The statement is obvious when $n=0,1$.
For $n=2$, we get
\begin{align}
\label{delta=2}
\begin{split}
\delta^{(2)}(xy)=&\delta^{(2)}(x)\delta^{(2)}(y)+
\delta^{(2)}(x)\big(\delta^{(1)}(y)^1+\delta^{(1)}(y)^2+\delta^{(0)}(y)^\emptyset\big)\\
&+\big(\delta^{(1)}(x)^1+\delta^{(1)}(y)^1+\delta^{(0)}(y)^\emptyset\big)\delta^{(2)}(y)\\
&+\delta^{(1)}(x)^1\delta^{(2)}(y)^2+\delta^{(1)}(x)^2\delta^{(2)}(y)^1,
\end{split}
\end{align}
so the statement also holds.

\noindent Assume that the statement is proved when $P$ is a $n$-tree.
Let $\bar{P}$ be a $(n+1)$-tree.
There exists an integer $k \in \{0,\dots,n-1\}$,
such that $\bar{P}$ may be viewed as the glueing of the $2$-tree on the $k$-th
leaf of a $n$-tree $P$.
Then we have
$$\delta^{(\bar{P})}=(\Id^{\otimes k}\otimes \delta^{(2)} \otimes \Id^{\otimes
n-k-1})
\circ \delta^{(P)}.$$
Let us assume, for instance, that $k=n-1$. If $\nu$ is an integer, set
$$S_\nu=\{(\Sigma',\Sigma'')|\Sigma',\Sigma''\subset \{1,\dots,\nu\}\hbox{ and
}\Sigma' \cup \Sigma''=\{1,\dots,\nu\}\}.$$
Then 
$$S_n=f_{\{n\},\emptyset}(S_{n-1}) \cup f_{\emptyset,\{n\}}(S_{n-1})\cup
f_{\{n\},\{n\}}(S_{n-1})\hbox{ (disjoint union)},$$
where $f_{\alpha,\beta}(\Sigma',\Sigma'')=(\Sigma'\cup\alpha,\Sigma''\cup
\beta)$. By hypothesis, we have
$$\delta^{(P)}(xy)=\sum\limits_{(\Sigma_1,\Sigma_2)\in
S_n}\delta^{(P_{\Sigma_1})}(x)^{\Sigma_1}
\delta^{(P_{\Sigma_2})}(y)^{\Sigma_2},$$
therefore
\begin{align*}
\delta^{(P)}(xy)=&\sum_{(\Sigma',\Sigma'')\in S_{n-1}}
\delta^{(P_{\Sigma'\cup\{n\}})}(x)^{\Sigma' \cup \{n\}}
\delta^{(P_{\Sigma''})}(y)^{\Sigma''}\\
&+\delta^{(P_{\Sigma'})}(x)^{\Sigma' }
\delta^{(P_{\Sigma''\cup\{n\}})}(y)^{\Sigma''\cup\{n\}}\\
&+\delta^{(P_{\Sigma'\cup\{n\}})}(x)^{\Sigma' \cup \{n\}}
\delta^{(P_{\Sigma''\cup\{n\}})}(y)^{\Sigma''\cup\{n\}}.
\end{align*}
Applying $\Id^{\otimes n-1} \otimes \delta^{(2)}$ to this identity and using
(\ref{delta=2}) and the identities
\begin{align*}
(\Id^{\otimes k}\otimes \delta^{(1)} \otimes \Id^{\otimes |P|-k-1})\circ
\delta^{(P)}&=\delta^{(P)},\\
(\Id^{\otimes k}\otimes \delta^{(0)} \otimes \Id^{\otimes |P|-k-1})\circ
\delta^{(P)}&=0,
\end{align*}
we get $\delta^{(\bar{P})}(xy)=$
\begin{align*}
\sum_{(\Sigma',\Sigma'')\in S_{n-1}}&
{ {\hbox{\huge (}}} \left( (\Id^{\otimes |\Sigma'|} \otimes \delta^{(2)}) \circ
\delta^{(P_{\Sigma' \cup \{n \}})}\right)(x)^{\Sigma' \cup \{n,n+1\}}
\delta^{(P_{\Sigma''})}(y)^{\Sigma''}\\
&\hskip-1.4cm+\delta^{(P_{\Sigma'})}(x)^{\Sigma'}
\left( (\Id^{\otimes |\Sigma''|} \otimes \delta^{(2)}) \circ
\delta^{(P_{\Sigma'' \cup \{n \}})}\right)(y)^{\Sigma'' \cup \{n,n+1\}}\\
&\hskip-1.4cm+\left( (\Id^{\otimes |\Sigma'|} \otimes \delta^{(2)}) \circ
\delta^{(P_{\Sigma' \cup \{n \}})}\right)(x)^{\Sigma' \cup \{n,n+1\}}
\left( (\Id^{\otimes |\Sigma''|} \otimes \delta^{(2)}) \circ
\delta^{(P_{\Sigma'' \cup \{n \}})}\right)(y)^{\Sigma'' \cup \{n,n+1\}}\\
&\hskip-1.4cm+\left( (\Id^{\otimes |\Sigma'|} \otimes \delta^{(2)}) \circ
\delta^{(P_{\Sigma' \cup \{n \}})}\right)(x)^{\Sigma' \cup \{n,n+1\}}
\left(
\delta^{(P_{\Sigma'' \cup \{n \}})}(y)^{\Sigma'' \cup \{n\}}+
\delta^{(P_{\Sigma'' \cup \{n \}})}(y)^{\Sigma'' \cup \{n+1\}}\right)\\
&\hskip-1.4cm+\left(
\delta^{(P_{\Sigma' \cup \{n \}})}(x)^{\Sigma' \cup \{n\}}+
\delta^{(P_{\Sigma' \cup \{n \}})}(x)^{\Sigma' \cup \{n+1\}}\right)
\left( (\Id^{\otimes |\Sigma''|} \otimes \delta^{(2)}) \circ
\delta^{(P_{\Sigma'' \cup \{n \}})}\right)(y)^{\Sigma'' \cup \{n,n+1\}}\\
&\hskip-1.4cm+\delta^{(P_{\Sigma' \cup \{n \}})}(x)^{\Sigma' \cup \{n\}}
\delta^{(P_{\Sigma'' \cup \{n \}})}(y)^{\Sigma'' \cup \{n+1\}}
+\delta^{(P_{\Sigma' \cup \{n \}})}(x)^{\Sigma' \cup \{n+1\}}
\delta^{(P_{\Sigma'' \cup \{n \}})}(y)^{\Sigma'' \cup \{n\}}{ {\hbox{\huge)}}}.
\end{align*}
So we get $\delta^{(\bar{P})}(xy)=$
\begin{align*}
\sum_{(\Sigma',\Sigma'')\in S_{n-1}}&
{ {\hbox{\huge(}}} \delta^{(\bar{P}_{\Sigma' \cup \{n ,n+1\}})}(x)^{\Sigma' \cup \{n,n+1\}}
\delta^{(\bar{P}_{\Sigma''})}(y)^{\Sigma''}\\
&\hskip-1.2cm+\delta^{(\bar{P}_{\Sigma'})}(x)^{\Sigma'}
\delta^{(\bar{P}_{\Sigma'' \cup \{n ,n+1\}})}(y)^{\Sigma'' \cup \{n,n+1\}}\\
&\hskip-1.2cm+\delta^{(\bar{P}_{\Sigma' \cup \{n ,n+1\}})}(x)^{\Sigma' \cup \{n,n+1\}}
\delta^{(\bar{P}_{\Sigma'' \cup \{n ,n+1\}})}(y)^{\Sigma'' \cup \{n,n+1\}}\\
&\hskip-1.2cm+\delta^{(\bar{P}_{\Sigma' \cup \{n ,n+1\}})}
(x)^{\Sigma' \cup \{n,n+1\}}
\left(
\delta^{(\bar{P}_{\Sigma'' \cup \{n\}})}(y)^{\Sigma'' \cup \{n\}}
+\delta^{(\bar{P}_{\Sigma'' \cup \{n+1\}})}(y)^{\Sigma'' \cup \{n+1\}}
\right)\\
&\hskip-1.2cm+\left(
\delta^{(\bar{P}_{\Sigma' \cup \{n\}})}(x)^{\Sigma' \cup \{n\}}
+\delta^{(\bar{P}_{\Sigma' \cup \{n+1\}})}(x)^{\Sigma' \cup \{n+1\}}
\right)
\delta^{(\bar{P}_{\Sigma'' \cup \{n ,n+1\}})}(y)^{\Sigma'' \cup \{n,n+1\}}\\
&\hskip-1.2cm+\delta^{(\bar{P}_{\Sigma' \cup \{n\}})}(x)^{\Sigma' \cup \{n\}}
\delta^{(\bar{P}_{\Sigma'' \cup \{n+1\}})}(y)^{\Sigma'' \cup \{n+1\}}
+\delta^{(\bar{P}_{\Sigma' \cup \{n+1\}})}(x)^{\Sigma' \cup \{n+1\}}
\delta^{(\bar{P}_{\Sigma'' \cup \{n\}})}(y)^{\Sigma'' \cup \{n\}}\!\!{ {\hbox{\huge)}}}.
\end{align*}
We have
\begin{align*}
S_{n+1}=&f_{\{n,n+1\},\{n,n+1\}}(S_{n-1})\cup f_{\{n,n+1\},\{n\}}(S_{n-1})
\cup f_{\{n,n+1\},\{n+1\}}(S_{n-1})\\
&\cup f_{\{n,n+1\},\emptyset}(S_{n-1})
\cup f_{\{n\},\{n,n+1\}}(S_{n-1})
\cup f_{\{n+1\},\{n,n+1\}}(S_{n-1})\\
&\cup f_{\emptyset,\{n,n+1\}}(S_{n-1})\cup f_{\{n\},\{n+1\}}(S_{n-1})
\cup f_{\{n+1\},\{n\}}(S_{n-1})\hbox{ (disjoint union)},\\
\end{align*}
where we recall that $f_{\alpha,\beta}(\Sigma',\Sigma'')=(\Sigma'\cup \alpha,\Sigma''\cup
\beta)$. So we get
$$\delta^{(\bar{P})}(xy)=\sum_{(\bar{\Sigma}',\bar{\Sigma}'')\in S_{n+1}}
\delta^{(P_{\bar{\Sigma}'})}(x)^{|\bar{\Sigma}'|}
\delta^{(P_{\bar{\Sigma}''})}(y)^{|\bar{\Sigma}''|}.$$
The proof is the same for a general $k \in \{0,\dots,n-1\}$.
This establishes the induction.
\findem

\vskip20pt

\centerline {\bf \S\; 6 \ Proofs of Proposition \ref{refinement}, 
Theorem \ref{theo:4} and Proposition \ref{metrized}}

\stepcounter{section}\label{section:6}

\vskip20pt

{\it 1. Proof of Proposition \ref{refinement}.} According to \cite{Dr2}, 
Proposition 3.10, there exists a series ${\cal E}'(\varphi)
\in U(\g)^{\otimes 3}[[\hbar]]$, expressed in terms of $(\mu,\varphi)$
by universal acyclic expressions (and therefore invariant), 
such that ${\cal E}'(\varphi) = 1 + O(\hbar^2)$, and ${\cal E'}(\varphi)$
satisfies the pentagon identity. Then $(U(\g)[[\hbar]],m_0,\Delta_0,
{\cal E}'(\varphi))$ is a quasi-Hopf algebra. By Theorem 
\ref{theo:3}, 2), 
there exists a twist $F \in U(\g)^{\otimes 2}[[\hbar]]^\times$, 
such that  $(U(\g)[[\hbar]],m_0,{}^F\Delta_0,
{}^F{\cal E}'(\varphi))$ is admissible.

${\cal E}(\varphi)$ gives rise to a collection of invariant 
elements ${\cal E}'(\varphi)_{p_1,p_2,p_3,n} \in \otimes_{i=1}^2 
S^{p_i}(\g)$, defined by the condition that the image of 
${\cal E}'(\varphi)$ by the symmetrization map 
$U(\g)^{\otimes 3}[[\hbar]] \to S^\cdot(\g)[[\hbar]]$ is 
$\sum_{n\geq 0,p_1,p_2,p_3\geq 0} \hbar^n 
{\cal E}'(\varphi)_{p_1,p_2,p_3,n}$. $F$ is then expressed using only 
the ${\cal E}'_{p_1,p_2,p_3,n}$, the Lie bracket and the symmetric group 
operations on the $\g^{\otimes n}$. So $F$ is invariant and defined by 
universal acyclic expressions. Therefore 
${}^F\Delta_0 = \Delta_0$. ${\cal E}(\varphi) := 
{}^F{\cal E}'(\varphi)$ is then expressed by universal 
acyclic expressions, and defines an admissible quantization of 
$(\g,\mu,\delta = 0,\varphi)$. 

\medskip 

{\it 2. Proof of Theorem \ref{theo:4}, 1).}
We have then ${\cal E}(\varphi) \in (U(\g)[[\hbar]]')^{\bar\otimes 3}$.  
Since the coproduct is $\Delta_0$, $U(\g)[[\hbar]]'$ is the complete 
subalgebra of $U(\g)[[\hbar]]$ generated by $\hbar \g$, so it is a 
flat deformation of $\widehat S^\cdot(\g)$ with Kostant-Kirillov Poisson 
structure. We then set $\widetilde\varphi := {\cal E}(\varphi)$ modulo 
$\hbar$.

\medskip 

{\it 3. Proof of Theorem \ref{theo:4}, 2).}
Let  $\widetilde{\varphi}_1$, $\widetilde{\varphi}_2$ be the elements of
$\widehat{S}^\cdot(\g)^{\bar\otimes 3}$ such that 
$$
(\widehat{S}^\cdot(\g),m_0,P_\g,\Delta_0,\widetilde{\varphi}_i)
$$
are Drinfeld algebras. Let $C$ be the lowest degree
component of  $\widetilde{\varphi}_1-
\widetilde{\varphi}_2$. Then the degree $k$ of $C$ is $\geq 4$.
Taking the degree $k$ part of the difference of the pentagon
identities for  $\widetilde{\varphi}_1$ and 
$\widetilde{\varphi}_2$, we find $\dif(C)=0$, where $\dif$~: 
$S^\cdot(\g)^{\otimes 3} \to S^\cdot (\g)^{\otimes 4}$ is
the co-Hochschild differential. So
$\Alt(C) \in \Lambda^3(\g)$, and
since $\Alt(C)$ also has degree $\geq 4$, $\Alt(C)=0$.
If $C_{p_1,p_2,p_3}$ is the component of $C$ in 
$\Otimes_{i=1}^{3} S^{p_i}(\g)$ then we may define inductively 
$B$ by $B_{0,k}=B_{1,k-1}=0$, $B_{2,k-2}={{1}\over{2}}
(\Id \otimes m)(C_{1,1,k-2})$, and
$$B_{i+1,k-i-1}={{1}\over{i+1}}(\Id \otimes m)[C_{i,1,k-i-1}+\big( (\Id \otimes
\dif)(B_{i,k-i})\big)_{i,1,k-i-1}],$$
where $B_{i,j}$ is the component of $B$ in $S^i(\g) \otimes S^j(\g)$
and $m$ is the product of $S^\cdot(\g)$.
So $B$ can be chosen to be $\g$-invariant.
Applying successive twists, we obtain the result.

\medskip
{\it 4. Proposition \ref{metrized}.} 
According to \cite{Dr3}, $(U(\g),m_0,\Delta_0,e^{\hbar t_\g/2},
\Phi(\hbar t_\g^{1,2},\hbar t_\g^{2,3}))$ is a quasi-triangular
quasi-Hopf algebra. One checks that it is admissible; then the
reduction modulo $\hbar$ of
the corresponding QFS algebra is the Drinfeld algebra of 1).
\findem

\medskip

\begin{remark}
In the proof of Theorem \ref{theo:4}, 2), 
we cannot use Theorem A of \cite{Dr2} because we do not
know that the twist constructed there is admissible.
\end{remark}

\vskip20pt

\centerline {\bf \S\; 7 \ Associators and Lie associators}

\stepcounter{section}\label{section:7}

\vskip20pt

In this section, we state precisely and prove Theorem \ref{theo:7.1}.

\bigskip

\noindent {\bf - a - Statement of the result}

\bigskip

Recall that the algebra $\TC_n$, $n \geq 2$, has generators $t^{i,j}$,
$1 \leq 1 \not=j \leq n$, and relations $t^{j,i}=t^{i,j}$,
\begin{align*}
[t^{i,j}+t^{i,k},t^{j,k}]&=0 \hbox{ when }i,j,k \hbox{ are all distinct},\\
[t^{i,j},t^{k,l}]&=0 \hbox{ when }i,j,k,l \hbox{ are all distinct}.
\end{align*}
$\tf_n$ is defined as the Lie algebra with the same generators and
relations. Then
$\TC_n=U(\tf_n)$.
When $n \leq m$ and $(I_1,\dots,I_n)$ is a collection of disjoint subsets of
$\{1,\dots,m\}$, there is a unique algebra morphism $\TC_n \to \TC_m$ taking
$t^{i,j}$ to
$\sum\limits_{\alpha \in I_i,\beta \in I_j}t^{\alpha,\beta}$. We call it an
insertion-coproduct morphism and denote it by $x \mapsto x ^{I_1,\dots,I_n}$.
In particular, we have an action of $\fgam_n$ on $\TC_n$. Let us attribute degree
$1$ to each generator $t^{i,j}$; this defines gradings on the algebra $\TC_n$
and on the Lie algebra $\tf_n$. We denote by ${\widehat{\TC}}_n$ and
${\widehat{\tf}}_n$ their completions for this grading.
Then ${\widehat{\TC}}_n$ is the preimage of $\KM^\times$ by the natural
projection
${\widehat{\TC}}_n\to \KM$, and the exponential
is a bijection
$({\widehat{\TC}}_n)_0 \to 1 + ({\widehat{\TC}}_n)_0$ (where
$({\widehat{\TC}}_n)_0=\Ker ({\widehat{\TC}}_n \to \KM)$). We have an exact
sequence
$$1 \to 1 + ({\widehat{\TC}}_n)_0 \to ({\widehat{\TC}}_n)^\times
\to \KM^\times \to 1.$$
An {\it associator} is an element $\Phi$ of $1 + ({\widehat{\TC}}_n)_0$,
satisfying the pentagon equation
\begin{equation}
\label{pent}
\Phi^{1,2,34}\Phi^{12,3,4}=\Phi^{2,3,4} \Phi^{1,23,4} \Phi^{1,2,3}, 
\end{equation}
the hexagon equations
\begin{equation*}
e^{{{t^{1,3}+t^{2,3}}\over{2}}}=\Phi^{3,1,2}e^{{{t^{1,3}}\over{2}}}
\big(\Phi^{1,3,2}\big)^{-1}e^{{{t^{2,3}}\over{2}}}\Phi^{1,2,3}
\end{equation*}
and
\begin{equation*}
e^{{{t^{1,2}+t^{1,3}}\over{2}}}=\big(\Phi^{2,3,1}\big)^{-1}
e^{{{t^{1,3}}\over{2}}}
\Phi^{2,1,3}e^{{{t^{1,3}}\over{2}}}\big(\Phi^{1,2,3}\big)^{-1}
\end{equation*}
and $\Alt(\Phi)={{1}\over{8}}[t^{1,2},t^{2,3}]+$ terms of degree $>2$.
We denote by $\Assocb$ the set of associators.
If $\Phi$ satisfies the duality condition $\Phi^{3,2,1}=\Phi^{-1}$,
then both hexagon equations are equivalent. We denote by $\Assocb^0$
the subset of all $\Phi \in \Assocb$ satisfying the duality
condition.
If $F\in 1 + ({\widehat{\TC}}_2)_0 $ and $\Phi \in 
1 + ({\widehat{\TC}}_3)_0 $, the {\it twist of $\Phi$ by $F$} is
$${}^F\!\Phi=F^{2,3}F^{1,23}\Phi(F^{1,2}F^{12,3})^{-1}.$$
This defines an action of $1 + ({\widehat{\TC}}_2)_0$ on 
$1 + ({\widehat{\TC}}_3)_0 $, which preserves
$\Pent=\{\Phi \in 1 + ({\widehat{\TC}}_3)_0|\Phi$ satisfies (\ref{pent})$\}$,
$\Assocb$ and $\Assocb^0$
($\Pent$ and $\Assocb$ are preserved because $F$ has the form $f(t^{1,2})$,
$f \in 1 + t\KM[[t]]$, so the ``twisted $R$-matrix'' ${}^F\! R=F^{2,1}R
F^{-1}=f(t^{2,1}) e^{t^{1,2}/2}$ $f(t^{1,2})^{-1}
=e^{t^{1,2}/2}$.
$\Assocb^0$ is preserved because each $F$ is such that $F=F^{2,1}$.)
We denote by $\Assocb_\Lie^0$, $\Assocb_\Lie$ and
$\Pent_\Lie$ the subsets
of all $\Phi$ in $\Assocb$, $\Assocb^0$ and $\Pent$, such that 
$\log(\Phi)\in {\widehat{\tf}}_3$.
\begin{theorem}
\label{thm:assoc}
There is exactly one element of $\Pent_\Lie$ resp., $\Assocb_\Lie$,
$\Assocb_\Lie^0$) in each orbit of the action of $1+({\widehat{\TC_2}})_0$
on $\Pent$ (resp., $\Assocb$, $\Assocb^0$).
The isotropy group of each element of $\Pent$ is
$\{e^{\lambda t^{1,2}}|\lambda \in \KM\}\subset 1+({\widehat{\TC_2}})_0$.

\end{theorem}

\bigskip

\noindent {\bf - b - Proof of Theorem \ref{thm:assoc}}

\bigskip

The arguments are the same in all three cases, so we treat the case of $\Assocb$.

\smallskip

Let $\Phi$ belongs to $\Assocb$. Set $\Phi=1 + \sum\limits_{i>0} \Phi_i$, where
$\Phi_i$ is the degree $i$ component of $\Phi$. 
Let $\dif$ be the co-Hochschild differential, 
\begin{align*}
\dif~:~\TC_n &\to \TC_{n+1}\\
x &\mapsto \sum\limits_{i=1}^{n}
{(-1)}^{i+1}x^{1,\dots,\{i,i+1\},\dots,n+1}-
x^{2,3,\dots,n+1}+{(-1)}^{n}x^{1,2,\dots,n}.
\end{align*}
Then $\dif(\Phi_2)=0$, and $\Alt(\Phi_2)={{1}\over{8}}[t^{1,2},t^{2,3}]$.
Computation shows that this implies that for some $\lambda \in \KM$, we
have $\Phi_2 ={{1}\over{8}}[t^{1,2},t^{2,3}]+
\lambda \dif((t^{1,2})^2)$.
We construct $F \in 1 + ({\widehat{\TC}})_0$, such that ${}^F\! \Phi
\in \Assocb_\Lie$, as an
infinite product $F=\cdots F_n \cdots F_2$, where $F_i \in 1 +
({\widehat{\TC}}_2)_{\geq i}$ (the index $\geq i$ means the part of degree $\geq
i$). If we set
$F_2=1 + \lambda (t^{1,2})^2$, then $\log({}^{F_2}\!\Phi) \in
\tf_3 + ({\widehat{\TC}}_3)_{\geq 3}$.
Assume that we have found $F_3,\dots,F_{n-1}$, such that
$\log({}^{{\bar{F}}_{n-1}}\! \Phi) \in \tf_3+
({\widehat{\TC}}_3)_{\geq n}$, where ${\bar{F}}_{n-1}=F_{n-1}\cdots F_2$.
Then $\varphi^{(n-1)}:=\log({}^{{\bar{F}}_{n-1}}\!\Phi)$ satisfies
$$\big( \varphi^{(n-1)}\big)^{1,2,34}\star 
\big(\varphi^{(n-1)} \big)^{12,3,4}=
\big(\varphi^{(n-1)} \big)^{2,3,4}\star
\big(\varphi^{(n-1)} \big)^{1,23,4}\star 
\big(\varphi^{(n-1)} \big)^{1,2,3},$$
where $\star$ is the CBH product in 
$({\widehat{\TC}}_3)_0$. Let $\varphi_n^{(n-1)}$ be the degree $n$ part of
$\varphi^{(n-1)}$. Then we get
$\dif(\varphi^{(n-1)})\in \tf_4$.
We now use the following satement, which will be proved in the next subsection.
\begin{proposition}
\label{prop:assoc}
If $\gamma \in \TC_3$ is such that $\dif(\gamma) \in \tf_4$, then there exists
$\beta \in \TC_2$, such that $\gamma + \dif(\beta) \in \tf_3$.
If $\gamma$ has degree $n$, one can choose $\beta$ of degree $n$.
\end{proposition}
\noindent
It follows that there exists $\beta \in \TC_2$ of degree $n$, such that
$\varphi_n^{(n-1)}- \dif(\beta) \in \tf_3$.
Set $F_n=1 + \beta$, then $\varphi^{(n)}=\log({}^{{\bar{F}}_{n}}\! \Phi)$
is such that 
$\varphi^{(n)}\in \varphi^{(n-1)}-\dif(\beta)+
({\widehat{\TC}}_3)_{\geq n+1}$, so
$\varphi^{(n)} \in \tf_3+ ({\widehat{\TC}}_3)_{\geq n+1}$.
Moreover, the product $F= \cdots F_n \cdots F_2$
is convergent, and 
${}^F\!\Phi$ then satisfies
$\log({}^F\!\Phi)\in \widehat{\tf}_3$. This proves the 
existence of $F$, such that 
${}^F\! \Phi \in \Assocb_\Lie$.

\smallskip

Let us now prove the unicity of an element of $\Assocb_\Lie$, twist-equivalent to
$\Phi \in \Assocb$. This follows from:
\begin{proposition}
\label{prop:unicity}
Let $\Phi'$ and $\Phi''$ be elements of $\Assocb_\Lie$, and let
$F$ belong to $1+({\widehat{\TC}}_2)_0$.
Then ${}^F\Phi'=\Phi''$ if and only if there exists $\lambda \in \KM$ such that
$F=e^{\lambda t^{1,2}}$ and $\Phi''=\Phi'$.
\end{proposition}
\noindent{\sc Proof of Proposition \ref{prop:unicity}.~}
Since $t^{1,2}+t^{1,3}+t^{2,3}$ is central in ${\widehat{\TC}}_3$, we have
${}^{F_\lambda }\!\Phi'=\Phi'$ when $F_\lambda=e^{\lambda t}$, for any $\lambda
\in \KM$.
Conversely, let $F_i$ be the degree $i$ part of $F$. Then for some $\lambda_0 \in \KM$, we
have
$F_1=\lambda_0 t$.
Replacing $F$ by $F'=FF_{-\lambda_0}$, we get
${}^{F'}\! \Phi' = \Phi''$,
and $F'-1$ has valuation $\geq 2$ (for the degree in $t$).
Assume that $F'-1 \not= 0$ and let $\nu$ be its valuation. Let $F_\nu'$ be the
degree $\nu$ part of $F'$.
Then 
$\dif(F_\nu')\in \tf_3$. On the other hand, $F_\nu'=\mu(t^{1,2})^\nu$,
where $\mu \in \KM-\{0\}$.
Now $\dif((t^{1,2})^\nu) \in \TC_3=U(\tf_3)$ has degree $\leq \nu$ for the
filtration
of $U(\tf_3)$, and its symbol in 
$S^\nu(\tf_3)=\gr_\nu(U(\tf_3))$ is
$\sum\limits_{\nu'=1}^{v-1}\begin{pmatrix}\nu \\ \nu' \end{pmatrix}
(t^{1,3})^{\nu '}(t^{2,3})^{\nu-\nu'}-
\sum\limits_{\nu''=1}^{\nu-1}\begin{pmatrix}
\nu \\ \nu'' \end{pmatrix}(t^{1,2})^{\nu''}(t^{1,3})^{\nu -\nu''}$: this is the image
of a non-zero element in $S^\nu(\KM t^{1,2} \oplus \KM t^{1,3} \oplus \KM
t^{2,3})$ under the injection $S^\nu (\Oplus_{1 \leq i<j \leq 3} \KM t^{i,j})
\hookrightarrow S^\nu(\tf_3)$, so it is non-zero. So $F'\not=1$ leads to a
contradiction. So $F=F_{\lambda_0}$, therefore
$\Phi''=\Phi'$.
\findem

\smallskip

\noindent
Note that we have proved the analogue of Proposition \ref{prop:assoc},
where the indices of $\TC_3,\tf_4$, etc., are shifted by $-1$.
\bigskip

\noindent {\bf - c - Decomposition of $\tf_3$ and proof of Proposition
\ref{prop:assoc}}

\bigskip

To end the proof of the first part of Theorem \ref{thm:assoc}, it remains to prove Proposition
\ref{prop:assoc}. For this, we construct a decomposition of $\tf_n$. For
$i=1,\dots,n$, there is a unique algebra morphism
$\epsilon_i$~: $ \TC_n \to \TC_{n-1}$,
taking $t_{i,j}$ to $0$ for any $j \not= i$, and taking $t_{j,k}$
to $t_{j-\lambda_i(j),k-\lambda_i(k)}$ if
$j,k \not=i$, where $\lambda_i(j)=0$ if $j < i$ and
$=1$ if $j >i$. Then $\epsilon_i$ induces a Lie algebra morphism
$\widetilde{\epsilon}_i$~: $\tf_n \to \tf_{n-1}$.
Set $\widetilde{\tf}_n=\CAP_{i=1}^n \Ker(\widetilde{\epsilon}_i)$.
Then we have

\begin{lemma}
$$\tf_n=\Oplus_{k=0}^n \hskip0.3cm \Oplus_{I \in \PC_k(\{1,\dots,
n\})}({\widetilde{\tf}}_k)^I,$$
where $\PC_k(\{1, \dots, n\})$ is the set of subsets of $\{1, \dots, n\}$ of
cardinal $k$, and
$({\widetilde{\tf}}_k)^I$ is the image of ${\widetilde{\tf}}_k$ under
$\tf_k \to \tf_n$, $x \mapsto x^{i_1,\dots,i_k}$,
where $I=\{i_1,\dots,i_k\}$.
\end{lemma}
\noindent{\sc Proof of Lemma.~}
Let $\Ff$ be the free Lie algebra with generators
$\widetilde{t}_{i,j}$, where $1 \leq i < j \leq n$.
It is graded by $\Gamma:=\NM^{\{(i,j)|1 \leq i <j \leq n \}}$: the
degree of $\widetilde{t}_{i,j}$ is the vector $\db_{i,j}$, whose $(i',j')$
coordinate
is $\delta_{(i,j),(i',j')}$. For
$\uk \in \Gamma$, we denote by $\Ff_{\uk}$ the part of $\Ff$ of degree $\uk$.
Let $\pi$~: $\Ff \to \tf_n$ be the canonical projection.
Since the defining ideal of $\tf_n$ is graded, we have
\begin{equation}
\label{tn:graded}
\tf_n=\Oplus_{\uk \in \Gamma} \pi (\Ff_{\uk}).
\end{equation}
On the other hand, one checks that 
$\widetilde{\tf}_n=\Oplus_{\uk \in \widetilde{\Gamma}} \pi(\Ff_{\uk})$, where
$\widetilde{\Gamma}$ is the set of maps $k$~: 
$\{(i,j)|1 \leq i<j \leq n\} \to \NM$, such that for each $i$,
$\sum\limits_{j|j>i} k (i,j)+ \sum\limits_{j|j<i}k(j,i) \not=0$.
Define a map $\lambda$~: $\Gamma \to \PC(\{1,\dots,n\})$ as
follows ($\PC(\{1,\dots,n\})$ is the set of subsets of $\{1,\dots,n\}$):
$\lambda$ takes the map $k$~: $\{(i,j)|1 \leq i<j \leq n\} \to \NM$
to 
$\{i|\sum\limits_{j|j>i} k(i,j) + \sum\limits_{j|j<i} k(j,i)\not=0\}$.
Then for each $I \in \PC(\{1,\dots,n\})$,
$\big(\widetilde{\tf}_{|I|}\big)^I$ identifies with 
$\Oplus_{\uk \in \lambda^{-1}(I)} \pi(\Ff_{\uk})$. Comparing with
(\ref{tn:graded}), we get
$$\hskip4.3cm\tf_n=\Oplus_{I \in \PC(\{1,\dots,n\})}\big(\widetilde{\tf}_{|I|}\big)^I.
\hskip4.3cm\findem$$
When $n=3$, we get $\tf_3=\KM t^{1,2} \oplus \KM t^{1,3} \oplus \KM t^{2,3}
\oplus
\widetilde{\tf}_3$.
On the other hand, the fact that the insertion-coproduct maps
take $\tf_n$ to $\tf_m$ implies that $\dif$~: $\TC_n \to \TC_{n+1}$ is compatible
with the filtrations induced by the identification 
$\TC_n=U(\tf_n)$, $\TC_{n+1}=U(\tf_{n+1})$. The associated graded map is
$$\gr^\cdot(\dif)~:~S^\cdot(\tf_n) \to S^\cdot(\tf_{n+1}).$$
Proposition \ref{prop:assoc} now follows from:
\begin{lemma}
When $k \geq 2$, the cohomology of the complex
$$S^k(\tf_2) \To^{\gr^k(\dif)} S^k(\tf_3) \To^{\gr^k(\dif)} S^k(\tf_4)$$
vanishes.
\end{lemma}
\noindent{\sc Proof of Lemma.~}
We have
\begin{equation}
\label{decomp}
S^k(\tf_3)=\Oplus_{\alpha=0}^kS^{k-\alpha}\big(\Oplus_{1 \leq i < j \leq 3} \KM
t^{i,j}\big)\otimes S^\alpha(\widetilde{\tf}_3).
\end{equation}
Let $x \in S^k(\tf_3)$, and let $(x_\alpha)_{\alpha=0,\dots,k}$ be its components
in the decomposition (\ref{decomp}).
We have
$$S^\cdot(\tf_4)=S^\cdot(\widetilde{\tf}_4)\otimes
\Bbo_{2 \leq i< j \leq 4} S^\cdot(\widetilde{\tf}_3^{1,i,j})
\otimes \Bbo_{i=2}^4 S^\cdot(\widetilde{\tf}_2^{1,i})\otimes
S^\cdot({\tf}_3^{2,3,4}).$$
We denote by $p$ the projection
$$p~:~S^\cdot({\tf}_4)\to \widetilde{\tf}_3^{1,3,4}\otimes
S^\cdot({\tf}_3^{2,3,4}),$$
which is the tensor product of: the identity on the last factor, the projection
to degree $1$ on the factor $S^\cdot(\widetilde{\tf}_3^{1,3,4})$, and the
projection
to degree $0$ in all other factors.
We also denote by $m$~: $\widetilde{\tf}_3^{1,3,4} \otimes
S^\cdot({\tf}_3^{2,3,4})\to S^\cdot(\tf_3)$ the map induced by the
identifications
$\widetilde{\tf}_3^{1,3,4} \subset {\tf}_3^{1,3,4}\simeq
\tf_3$, ${\tf}_3^{2,3,4}\simeq
\tf_3$ followed by the product map in $S^\cdot(\tf_3)$.
We denote by $\dif_1,\dif_2,\dif_3$ the maps $\TC_3 \to \TC_4$ defined by
\begin{align*}
\dif_1(x)&=x^{12,3,4}-x^{1,3,4}-x^{2,3,4},\\
\dif_2(x)&=x^{1,23,4}-x^{1,2,4}-x^{1,3,4},\\
\dif_3(x)&=x^{1,2,34}-x^{1,2,3}-x^{1,2,4},
\end{align*}
so $\dif=\dif_1-\dif_2+\dif_3$. The maps $\dif_i$ are compatible with the
filtrations of $\TC_3$ and $\TC_4$; we denote by $\gr^k(\dif_i)$ the
corresponding graded maps, so $\gr^k(\dif)=\gr^k(\dif_1)-
\gr^k(\dif_2)+\gr^k(\dif_3)$.
Then if we set 
$$x_1=\sum_{a,b,c|a+b+c=k-1}(t^{1,2})^a(t^{1,3})^b(t^{2,3})^c \otimes
e_{a,b,c},$$
where $e_{a,b,c} \in \widetilde{\tf}_3$, we have
$$
m \circ p \circ \gr^k(\dif_1)(x)= \big( \sum_{\alpha=0}^k \alpha x_\alpha \big) 
-(t^{2,3})^{k-1} e_{0,0,k-1}.
$$
On the other hand, let us define the $i$-degree of an element of
$\big(\widetilde{\tf}_{|I|}\big)^I$ to be $1$ if $i\in I$ and $0$ if $i \notin I$. Then
the $i$-degree of $\Otimes_{I \subset \{1,\dots,n\}} S^{\alpha_I}
\big(\big(\widetilde{\tf}_{|I|}\big)^I\big) \subset S^\cdot(\tf_n)$ is $\sum\limits_{I|i\in
I}\alpha_I$. If $x$ is homogeneous for the $1$-degree, then so is
$\gr^k(\dif_2)(x)$, and $1$-degree$(\gr^k(\dif_2)(x))=1$-degree$(x)$.
On the other hand, the elements of $S^\cdot(\tf_4)$ whose $1$-degree is $\not=1$
are in the kernel of $p$.
It follows that
$$m \circ p \circ \gr^k(\dif_2)(x_\alpha)=0 \hbox{ if } \alpha \not=1,$$
and
$p \circ \gr^k(\dif_2)(x_1)=(e_{0,0,k-1})^{1,3,4}\big[
(t^{2,4}+t^{3,4})^{k-1}-(t^{3,4})^{k-1}\big],$
so 
$$m\circ p \circ \gr^k(\dif_2)(x_1)=e_{0,0,k-1}\big[(t^{1,3}+t^{2,3})^{k-1}
-(t^{2,3})^{k-1}\big].$$ 
Finally, $p \circ \gr^k(\dif_3)(x)=0.$
If $x$ is such that 
$\gr^k(\dif)(x)=0$, we have $m \circ p \circ \gr^k(\dif)(x)=0$,
so
$$\sum_{\alpha \geq 0} \alpha
x_\alpha=e_{0,0,k-1}(t^{1,3}+t^{2,3})^{k-1}.$$
Looking at degrees in the decomposition (\ref{decomp}), we get $x_\alpha=0$ for
$\alpha \geq 2$, and
$x_1=e_{0,0,k-1}(t^{1,3}+t^{2,3})^{k-1}$.
Using the projection
$p'$~: $S^\cdot(\tf_4)\to \widetilde{\tf}_3^{1,2,4} \otimes
S^\cdot(\tf_3^{1,2,3})$,
we get in the same way
$x_1=e_{k-1,0,0}(t^{1,2}+ t^{1,3})^{k-1}$.
Now $e_{k-1,0,0}(t^{1,2}+t^{1,3})^{k-1}=
e_{0,0,k-1}(t^{1,3}+t^{2,3})^{k-1}$ implies
$e_{k-1,0,0}=e_{0,0,k-1}=0$ so $x_1=0$.
Therefore $x \in S^k\big( \Oplus_{1 \leq i < j \leq
3} \KM t^{i,j}\big)$. Let us set
$x=S(t^{1,2},t^{1,3},t^{2,3})$, where $S$ is a homogeneous
polynomial
of degree $k$ of $\KM[u,v,w]$.
Since $\dif(x)=0$, we have
\begin{multline*}
S(t^{1,3}+t^{2,3},t^{1,4}+t^{2,4},t^{3,4})-
S(t^{1,2}+t^{1,3},t^{1,4},t^{2,4}+t^{3,4})\\
+S(t^{1,2},t^{1,3}+t^{1,4},t^{2,3}+t^{2,4})
=S(t^{2,3},t^{2,4},t^{3,4})+S(t^{1,2},t^{1,3},t^{2,3})
\end{multline*}
(equality in $S^\cdot \big( \Oplus_{1 \leq i < j \leq
4} \KM t^{i,j}\big)$).

Applying ${{\partial}\over{\partial t^{1,2}}} \circ 
{{\partial}\over{\partial t^{3,4}}}$ to this equality, we get
$$(\partial_u\partial_w S)(t^{1,2}+t^{1,3},t^{1,4},t^{2,4}+t^{3,4})=0,$$
therefore $\partial_u \partial_w S=0$. We have therefore
$$S(u,v,w)=P(u,v)+Q(v,w),$$
where
$P$ and $Q$ are homogeneous polynomials of degree $k$.
Moreover, $\dif(x)=0$, so
\begin{align}
\label{final:id}
\begin{split}
&\big[P(t^{1,2},t^{1,3}+t^{1,4})-P(t^{1,2}+
t^{1,3},t^{1,4})-P(t^{1,2},t^{1,3})\big]\\
&+\big[Q(t^{1,4}+t^{2,4},t^{3,4})-Q(t^{1,4},t^{2,4}
+t^{3,4})-Q(t^{2,4},t^{3,4})\big]\\
&+\big[P(t^{1,3}+t^{2,3},t^{1,4}+t^{2,4})+Q(t^{1,3}+t^{1,4},t^{2,3}+t^{2,4})-
P(t^{2,4},t^{2,4})-Q(t^{1,3},t^{2,3})\big]=0.
\end{split}
\end{align}
Write this as an identity
$$B(t^{1,2},t^{1,3},t^{1,4})
+C(t^{1,4},t^{2,4},t^{3,4})+A(t^{2,3},t^{1,4},t^{1,3},t^{2,4})=0.$$
Then $A$ (resp., $B,C$) is independent on $t^{2,3}$ (resp., $t^{1,2}$,
$t^{3,4}$). Let us now determine $P$ and $Q$.
Since $B(t^{1,2},t^{1,3},t^{1,4})=B(0,t^{1,3},t^{1,4})$, we have
$P(u,v+w)-P(u+v,w)-P(u,v)=P(0,v+w)-P(v,w)-P(0,v).$
Therefore $(\dif \widetilde{P})(u,v,w)=0$, where
$\widetilde{P}(u,v)=P(u,v)-P(0,v)$ and
$\dif$ is the co-Hochschild differential of polynomials in one variable.
The corresponding cohomology is zero, so we have a polynomial
$\bar{P}$, such that
$$P(u,v)-P(0,v)=\bar{P}(u+v)-\bar{P}(u)-\bar{P}(v).$$
We conclude that $P(u,v)$ has the form
\begin{equation}
\label{form:P}
P(u,v)=\bar{P}(u+v)-\bar{P}(u)-R(v)
\end{equation}
where $\bar{P}$ and $R$ are polynomials in one variable of degree $k$; since
$P(u,v)$ is homogeneous of degree $k$, we can assume that $\bar{P}$ and $R$ are
monomials of degree $k$.
In the same way, since $C(t^{1,4},t^{2,4},t^{3,4})=C(t^{1,4},t^{2,4},0),$
we have
$Q(u+v,w)-Q(u,v+w)-Q(v,w)=Q(u+v,0)-Q(u,v)-Q(v,0),$
so $(\dif\widetilde{Q})(u,v,w)=0$, where $\widetilde{Q}(u,v)=Q(u,v)-Q(u,0)$.
So $Q(u,v)$ has the form
\begin{equation}
\label{form:Q}
Q(u,v)=\bar{Q}(u+v)-\bar{Q}(v)-S(u),
\end{equation}
where $\bar{Q}$ and $S$ are polynomials in one variable of degree $k$, which can 
be assumed to be monomials of degree $k$.
We have therefore
$$x=\bar{P}^{1,23}+\bar{Q}^{12,3}-\bar{P}^{1,2}-\bar{Q}^{2,3}-T^{1,3},$$
where $\bar{P}=\bar{P}(t^{1,2})$, $\bar{Q}=\bar{Q}(t^{1,2})$ and
$T=(R+S)(t^{1,2}).$ So 
$x=\dif(\bar{Q})+(\bar{P}+\bar{Q})^{1,23}-(\bar{P}+\bar{Q})^{1,2}-T^{1,3}.$
Set $a=\bar{P}+\bar{Q}$; we have $\dif(y)=0$, where
$y=a^{1,23}-a^{1,2}-T^{1,3}$; applying $\epsilon_1$ to $\dif(y)=0$,
we get
$T^{2,3}-T^{2,4}=0$, so $T=0$. We then get
$a^{12,34}-a^{12,3}-a^{2,34}+a^{2,3}=0$.
Applying $\epsilon_3 \circ \epsilon_2$ to this identity, we get
$a^{1,4}=0$.
Finally $\bar{P}=-\bar{Q}$, so $x =\dif(\bar{Q})$, which proves the lemma.
\findem

\bigskip

\noindent {\bf - c - Isotropy groups}

\bigskip

Proposition \ref{prop:unicity} can be generalized to the case of a pair of
elements of $\Pent_\Lie$, and it implies that the isotropy
group of each element of $\Pent_\Lie$ is
the additive group
$\{e^{\lambda t^{1,2}},\lambda \in \KM\}$. Let $\Phi$ be an element of
$\Pent$. There exists an element $\Phi_\Lie$ of $\Pent_\Lie$ in the
orbit
of $\Phi$. So the isotropy groups of
$\Phi$ and $\Phi_\Lie$ are conjugated. Since
$1+({\widehat{\TC}}_2)_0$ is commutative, the isotropy group of $\Phi$ is
$\{e^{\lambda t^{1,2}},\lambda \in \KM\}$.
\findem

\vskip1.1truecm

\centerline { ACKNOWLEDGEMENTS }

\vskip18pt

B.E. would like to thank Y. Kosmann-Schwarzbach for communicating to him
\cite{Dr4}, where the question of the construction of quasi-groups is raised, as well 
as for discussions on this question in 1995.
Both authors would also like to thank P. Etingof, F. Gavarini, Y. Kosmann-Schwarzbach
and P. Xu
for discussions.

\end{document}